\documentclass[11pt]{article}
\usepackage{amssymb,amsmath,latexsym}
\usepackage{epsfig}
\usepackage{enumerate}

\usepackage{color}

\renewcommand{\baselinestretch}{1.1}

\def\nn{{\nonumber}}
\newtheorem{thm}{Theorem}[section]

\newtheorem{lemma}[thm]{Lemma}
\newtheorem{prop}[thm]{Proposition}

\newtheorem{assumption}[thm]{Assumption}
\newtheorem{remark}[thm]{Remark}
\newtheorem{example}[thm]{Example}

\numberwithin{equation}{section}

\def\pf{{\medskip\noindent {\bf Proof. }}}
\def\qed{{\hfill $\Box$ \bigskip}}

\def\sD {{\cal D}} \def\sE {{\cal E}} \def\sF {{\cal F}}
  
  \def\sL {{\cal L}}
 \def\sN {{\cal N}}

 \def\bB {{\mathbb B}} 
 \def\bE {{\mathbb E}}

\def\bP {{\mathbb P}}  \def\bR {{\mathbb R}}

\def\wh{\widehat}
\def\<{\langle}
\def\>{\rangle}

\def\E{{\mathbb E}}
\def\P{{\mathbb P}}

\def\bea{\begin{align*}}
\def\eea{\end{align*}}
\def\bee{\begin{equation}}
\def\eee{\end{equation}}

\def\wh{\widehat}

\makeatletter
\@addtoreset{equation}{section}

\makeatother

\begin{document}
\bibliographystyle{plain}

\title{\Large \bf
Time Fractional Poisson Equations: Representations and Estimates}
\author{Zhen-Qing Chen \quad Panki Kim \quad Takashi Kumagai
\quad Jian Wang}

\date{}

\maketitle

\begin{abstract}
In this paper, we study existence and uniqueness of
strong as well as weak
solutions for general time fractional Poisson equations. We show that there is an
integral representation of the solutions of time fractional Poisson equations with zero initial values
in terms of semigroup for the  infinitesimal spatial generator  $\sL$
and the corresponding subordinator associated with the time fractional derivative.
This integral representation has an integral kernel $q(t, x, y)$, which we call
the fundamental solution for the time fractional Poisson equation, if the semigroup
for $\sL$ has an integral kernel. We further show that $q(t, x, y)$
can be expressed as a time fractional derivative
of the fundamental solution for the homogenous time fractional equation
under the assumption that the associated subordinator admits a conjugate subordinator.
Moreover,
when the Laplace exponent of
the associated subordinator satisfies the weak scaling property and its distribution  is self-decomposable,
we establish two-sided estimates for
the fundamental solution $q(t,x, y)$
through explicit estimates of
transition density functions of subordinators.
\end{abstract}

\bigskip
\noindent {\bf AMS 2010 Mathematics Subject Classification}: Primary
26A33, 60H30; Secondary:    34K37

\bigskip\noindent
{\bf Keywords and phrases}: time fractional derivative, Caputo derivative, Poisson equation,
subordinator, inverse subordinator, conjugate subordinator,
 fundamental solution

\bigskip
\allowdisplaybreaks

\section{Introduction}
In recent years, linear as well as non-linear
 partial differential equations with
fractional derivatives are extensively studied.
 This interest has been stimulated by
anomalous diffusion phenomena appeared in diverse
 fields including mathematics, physics, chemistry,
engineering, biology, geophysics and hydrology,
see \cite{MSi, MR, SKM, Uvv1, Uvv2} and the references therein.
The present paper is
concerned with solutions of general time fractional Poisson equations.

For a given function
$w: (0,\infty)\to [0,\infty)$ that
is unbounded, non-increasing and having $\int_0^\infty \min\{1, x\} (-d w(x))<\infty$,
a generalized time fractional  derivative
with weight $w$ is defined by
\begin{align}
\label{e:derfw}
\partial^w_t f(t):= \frac{d}{dt} \int_0^t w(t-s) (f(s)-f(0))\, ds,
\end{align}
 whenever the right hand side is well defined.
 See \cite{C, CKKW, HeK, K}.
 In particular,
  when $w(s)=\frac{1}{\Gamma(1-\beta)} s^{-\beta}$ for $\beta\in(0,1)$ (where $\Gamma(t)=\int_0^\infty s^{t-1}e^{-s}ds$ is the Gamma function),
 $\partial^w_t f$ is just the Caputo derivative of order $\beta$, i.e.,
 $$\partial^\beta_t f(t):= \frac{1}{\Gamma(1-\beta)}\frac{d}{dt} \int_0^t (t-s)^{-\beta} (f(s)-f(0))\, ds.$$

Let $\sL$ be a linear (can be unbounded)
operator on some function space over a locally compact separable metric space $E$ with a fully-supported Radon measure $\mu$.
We consider the following general time fractional Poisson equation
\begin{equation}\label{e:3.2}
 \partial^w_t u (\cdot , x) (t)= \sL u (t, \cdot) (x)+ f(t, x)
\quad t\in (0,  T_0], \ x \in E.
    \end{equation}
where $T_0>0$ and  $f(t,x)$ is a given function
on $(0,  T_0] \times E$.
For simplicity, in the following we will write the equation in \eqref{e:3.2} as
$$
 \partial^w_t u(t, x) = \sL u(t, x) + f(t, x) .
$$
When
$f(t,x)=0$, \eqref{e:3.2} is reduced into the following general time fractional Cauchy equation
\begin{equation}\label{e:3.2-0}
    \partial^w_t u(t, x) = \sL u(t, x),\quad
t\in (0,  T_0], \ x \in E.
    \end{equation}
Solutions of time fractional
Cauchy as well as Poisson
 equations have been attracted a lot of
attentions in the community of analysis, PDEs and stochastic analysis, see 
\cite{ACV, BM, CKK0, EIK, FN, GJ, LPJ, MN, U,US}
and the references therein. We note that most of
 quoted papers are concentrated on the Caputo derivative of fractional order.
Recently, Chen \cite{C} studied equation \eqref{e:3.2-0} with $T_0=\infty$ for any unbounded
and non-increasing function $w$ with $\int_0^\infty  \min\{1,x\}
\, (-dw(x)) <\infty$
and for any infinitesimal generator $\sL$
that generates a
 uniformly bounded and strongly continuous
semigroup in some Banach space,
and obtained strong existence and uniqueness
of solutions to \eqref{e:3.2-0}. The main feature of the approach in \cite{C} is
 a detailed analysis of the subordinator associated with the function $w$ together
with a probabilistic representation of the solution.
The existence and uniqueness results in \cite{C} for strong solutions of  \eqref{e:3.2-0} have been extended
to that of weak solutions in \cite{CKKW} when the infinitesimal generator
$\sL$ is a symmetric operator in some $L^2$ space
with $T_0=\infty$.
Moreover, two-sided estimates on the fundamental solution of \eqref{e:3.2-0} have been derived in \cite{CKKW} for a
family of $w$ and $\sL$.

In this paper, we call $q(t,x,y)$ the fundamental solution to the time fractional Poisson equation
\eqref{e:3.2} (with $u(0, x)=0$),  and  call $p(t,x,y)$ the fundamental solution to the homogenous
time fractional equation \eqref{e:3.2-0}.

The aim of our paper is twofold.
One is to study existence and uniqueness of
strong as well as weak  solutions for \eqref{e:3.2} and
to derive an integral representation of the solution to
the time fractional Poisson equation \eqref{e:3.2}
with zero initial data in terms of the semigroup $\{P^0_t; t\geq 0\}$
of $\sL$ and the density $\bar p(r, t)$ of the   driftless subordinator $S_r$
  whose L\'evy measure is $\nu$
 and $   w(x)=\nu (x, \infty)$ is the non-increasing function $w$ in \eqref{e:derfw}.
When $\{P^0_t; t\geq 0\}$ has an integral kernel $p^0 (t, x, y)$,
the integral representation obtained in this paper
for the solution of \eqref{e:3.2} admits an integral kernel $q(t, x, y)$,
which we call the the fundamental solution to the time fractional Poisson equation
\eqref{e:3.2} (with $u(0, x)=0$).
The other is to establish two-sided estimates on $q(t, x, y)$.
Our approach is mainly based on probabilistic ideas combined with some analytic techniques.
The paper can be viewed as a sequel to \cite{C, CKKW}, where existence and uniqueness,
probabilistic representations as well as two-sided estimates for
the fundamental solution to the homogenous time fractional equation
 \eqref{e:3.2-0}
are given.
In particular, it is shown in \cite{C} that $p(t, x, y):=\E \left[ p^0(E_t, x, y) \right]$,
where $E_t$ is the inverse of the subordinator $S$,  is the fundamental solution
for the time-homogenous time fractional equation \eqref{e:3.2}. However,
the study of the  fundamental solution $q(t, x y)$ of  \eqref{e:3.2} is  harder
than the fundamental solution $p(t, x, y)$
  of the equation \eqref{e:3.2-0}.
To obtain two-sided estimates for
   the fundamental solution $q(t, x, y)$ for the time fractional Poisson equation, we also need
some explicit estimates for transition density functions of
   subordinators,
 which are highly non-trivial.

Two recent papers \cite{HKT, To} considered a slightly different form of time fractional
derivative but  have a probabilistic representation
analogous to us. The proofs there are different from ours and scope of that two
 papers are restricted to
  Feller generators in Euclidean spaces.
	See also \cite{DTZ} for a recently related work
	where the time fractional derivative having a possibly time-dependent
	kernel and the spatial
	operator is a Dirichlet Laplacian in a regular Euclidean domain.

\medskip

We now present the precise setting of this paper.
Throughout this paper, let $S=\{S_t, \P; t\geq 0\}$ be a driftless
subordinator with $S_0=0$ 
that has a 
bounded density $\bar p(t, \cdot)$ for every $t>0$.
Denote by  {$\phi$} the {Laplace  exponent of $S$}; that is,
$$
\E \left[ e^{-\lambda S_t} \right] = e^{-t \phi (\lambda)}, \quad \lambda >0, t\ge0.
$$
The Laplace  exponent $\phi$ of $S$ is also called Bernstein function (vanishes at the origin) in the literature.
Since $S$ has no drift, it is well known
that there is a unique
measure $\nu$  on $(0, \infty)$, which is called L\'evy measure of $S$, satisfying
 $\int_0^\infty
   \min\{1,x\} \,\nu (dx)<\infty$
such that
$$ \phi (\lambda) = \int_0^\infty (1-e^{-\lambda x})\, \nu ( dx).
$$
Let
 \begin{align}
 \label{e:defw}
  w(x):=\nu (x, \infty) .
 \end{align}
We always assume that the L\'evy measure $\nu$ of the subordinator $S$ is infinite, which is equivalent to
$w(x)$ being unbounded.
  (However, $w$ is always
	locally integrable on $[0,\infty)$; see
  \cite[(2.2)]{C}.)
 Let $(\sL, \sD (\sL))$ be an infinitesimal generator corresponding to
a   uniformly bounded and strongly continuous
 semigroup $\{P_t^0;t\ge 0\}$ in some
  Banach   space. In this paper, we are concerned with  time fractional Poisson equation \eqref{e:3.2} with the weight function $w(x)$ and the infinitesimal generator $(\sL, \sD (\sL))$.

We now state
the existence and uniqueness for
the strong as well as weak solutions of
time fractional heat equations, which is one of the main results of this paper.
Define for $t>0$,
$E_t =\inf\left\{ s>0: S_s >t\right\}$, the inverse subordinator.

\begin{thm}[{\bf Strong solution}]\label{T:mains}
Assume that
the subordinator $\{S_t, \P; t\geq 0\}$
      has a bounded density $\bar p(t, \cdot )$ for each $t>0$
  and its L\'evy measure $\nu$ is infinite.
Suppose that $\{P^0_t; \, t\geq 0\}$ is a
 uniformly bounded and strongly
continuous semigroup in some Banach space $(\bB, \| \cdot \|)$ over a locally compact separable metric space $E$ and
$(\sL, \sD (\sL))$ is its infinitesimal generator. Let
$T_0\in(0,\infty)$, and
 $g\in \sD (\sL)$  and $f(t, x)$ be
a function defined on
$ (0, T_0] \times E$ so that for a.e. $t\in (0, T_0]$,
$f(t, \cdot) \in \sD (\sL)$ and
\begin{equation}\label{e:1.5}
\int_0^{T_0}  \| \sL f(t, \cdot) \| \, dt <\infty
\quad \hbox{ and } \quad \| f(t, \cdot )\| \leq K <\infty \hbox{ for a.e. } t\in (0,  T_0].
\end{equation}
Then the function
\begin{equation} \label{e:defu}\begin{split}
u(t,x):&= \bE \left[ P^0_{E_t} g (x) \right]+ \E \left[\int_{r=0}^\infty
{\bf1}_{\{ S_r<t\}} P^0_r f (t-S_r, \cdot) (x)\, dr \right]  \\
&= \bE \left[ P^0_{E_t} g (x) \right]+ \int_{s=0}^t \int_{r=0}^\infty P^0_r f(s,
\cdot)(x)\bar p(r, t-s)\,dr \, ds
\end{split}\end{equation}
 is the unique
$($strong$)$ solution of
   $\partial^w_t u (t, x) = \sL u (t, x)+ f(t, x)$
  on $(0,
   T_0]\times E$ with $u(0, x)=g(x)$
in the following sense:

\begin{description}
\item{\rm (i)} \,\, $u(t, \cdot)$ is well defined as an element in $\bB$ for each
$t\in (0,  T_0]$
such that $\sup_{t\in (0,  T_0]} \| u(t, \cdot)\|
<\infty$, $t\mapsto u(t,x)$ is continuous in $(\bB,\|\cdot\|)$ and $\lim_{t\to 0} \| u(t, \cdot )-g(\cdot)\| =0$.

\item{\rm (ii)} For a.e. $t\in (0,  T_0]$, $u(t, \cdot) \in \sD (\sL)$ and  $\sL u(t, \cdot)$ exists in the Banach space $\bB$
such that $\int_0^{T_0} \| \sL u(t, \cdot) \| \,dt <\infty$.

\item{\rm (iii)}
For every $T \in (0,   T_0]$, $\int_0^T w(T-t) (u(t, x) -g(x)) \,dt$ is absolutely convergent in
$(\bB, \| \cdot \|)$ and
$$
\int_0^T w(T-t) \left(u(t, \cdot )-g(\cdot) \right)\,dt = \int_0^T \left( f(t, \cdot) +
\sL u(t, \cdot )\right) \,dt \quad \hbox{in } \bB.
$$
\end{description}
\end{thm}

\medskip

\begin{remark} \rm
\begin{description}
\item{(i)} Clearly by Theorem \ref{T:mains},
if condition \eqref{e:1.5} holds for every $T_0>0$, then the function $u(t, x)$ given
by \eqref{e:defu} is well defined for all
$t >0$
and is the unique strong solution of
 $\partial^w_t u (t, x) = \sL u (t, x)+ f(t, x)$
  on $(0,\infty) \times E$ with $u(0, x)=g(x)$. Similar remark applies to Theorem \ref{T:mainw} below
	on weak solutions as well.

\item{(ii)}
When   $f=0$, the above result is established in \cite[Theorem 1.1]{C}.  We point out that in
\cite[Theorem 2.3]{C} (where $f=0$), the continuity of $t\mapsto \sL u(t, \cdot)$  in $(\bB, \| \cdot \|)$ appears as a part of the definition of the strong  solution. But one can see from  the uniqueness part of  the proof of \cite[Theorem 2.3]{C} that the continuity of $t\mapsto \sL u(t, \cdot)$  in $(\bB, \| \cdot \|)$  is not used.
Thus the continuity of $t\mapsto \sL u(t, \cdot)$
 should be one of properties of the strong solution
there, and can be removed from the definition of
 the strong solution in the statement of \cite[Theorem 2.3]{C}.
\end{description}
\end{remark}

\medskip

Note that, if we extend the definition of $f(t, \cdot)$ to $(-\infty,  T_0]$ by letting $f(t, \cdot)\equiv 0$ for $t \le 0$, then
\eqref{e:defu} can be rewritten as
\begin{align}
\label{e:defu0}
u(t,x)= \bE \left[ P^0_{E_t} g (x) \right]+ \bE \left[ \int_{r=0}^\infty P^0_r f(t-S_r,
\cdot)(x) \, dr \right].
\end{align}
(Cf., \cite[(35)]{HKT} and  \cite[(1.3)]{To}.)

\begin{thm} [{\bf Weak solution}]\label{T:mainw}
Assume that
the subordinator $\{S_t, \P; t\geq 0\}$
  has a bounded density $\bar p(t, \cdot )$ for each $t>0$ and its L\'evy measure $\nu$ is infinite.
 Suppose that $E$ is a locally compact separable metric
  space  and that $\mu$ is a $\sigma$-finite Radon measure on $E$ with full support.
 Suppose that
  $(\sL, \sD (\sL))$ is the
infinitesimal generator of a Dirichlet form on $L^2(E; \mu)$
and
  $\{P^0_t; t\geq 0\}$ is its associated transition
 semigroup in $L^2(E; \mu)$.
 Denote by $\langle \cdot,\cdot\rangle$ the inner product of
$L^2(E;\mu)$.
Suppose that
$T_0\in (0,\infty)$,
and $g \in  {L^2(E; \mu)}$ and that $f(t, x)$ is a function on $(0,  T_0]\times E$ so that
$$\| f(t, \cdot )\|_{L^2(E; \mu)} \leq K<\infty \ \hbox{ for a.e. } t\in (0,  T_0].$$
Then $u(t,x)$ defined by \eqref{e:defu} is the unique weak solution of
 $\partial^w_t u(t, x) = \sL u(t, x) + f(t, x)$
on $(0,  T_0]\times E$ with $u(0,
x)=g(x) $ in the following sense:

\begin{description}
\item{\rm (i)} $u(t, \cdot)$ is well defined as an element in $L^2(E; \mu)$ for every $t\in (0,  T_0]$
satisfying that $t\mapsto u(t,x)$ is continuous in
$L^2(E; \mu)$,
$$\sup_{t\in (0,  T_0]} \| u(t, \cdot )\|_{L^2(E; \mu)}  <\infty
\quad \hbox{and} \quad u(t, \cdot) \to
g \hbox{ in } L^2(E; \mu) \hbox{  as  } t\to 0.
$$

\item{\rm (ii)}
For every
$T\in (0,  T_0]$,
 $\int_0^T w(T-t) (u(t, x) -g(x)) \,dt$ is absolutely convergent in $L^2(E; \mu)$, and
\begin{equation}\label{e:Tmainw}
\left\< \int_0^T w(T-s) \left( u(s, \cdot) -g (\cdot) \right)\,ds, \varphi \right\>
=\int_0^T\<f(s, \cdot), \varphi \> \,ds + \int_0^T \< u(s, \cdot), \sL
\varphi \> \,ds
\end{equation}
for every $\varphi \in \sD (\sL)$.
\end{description}
Furthermore, if   $\sup_{t\in (0,  T_0]} \| \sL f(t, \cdot)\|
<\infty$, then $t\mapsto \sL u(t,x)$ is also  continuous in $L^2(E; \mu)$.
\end{thm}

A direct but important consequence of Theorems \ref{T:mains} and \ref{T:mainw} is that,
when the semigroup $\{P^0_t; t\geq 0\}$ has a density function $p^0(t, x, y)$,
 the fundamental solution $q(t, x, y)$ for the time fractional Poisson equation is
  \begin{align}\label{qp0p}
    q(t, x, y):= \int_0^\infty p^0(r, x, y) \bar p(r, t) \,dr.
  \end{align}
Formula
\eqref{qp0p} clearly implies the positivity of $q(t, x, y)$.
Furthermore, \eqref{qp0p} enables us to establish in Section \ref{section4}
two-sided estimates for
 the fundamental solution $q(t, x, y)$ for the time fractional Poisson equation
using estimates of $p^0(r, x, y)$ and $\bar p(r, t)$ only.

 Throughout the paper,
we write
$h(s)\simeq f(s)$,
if there exist constants $c_{1},c_{2}>0$ such that
$
c_{1}f(s)\leq h(s)\leq c_{2}f(s)
$
for the specified range of the argument $s$. Similarly, we write
$h(s)\asymp f(s)g(s)$,
if there exist constants $C_{1},c_{1},C_{2},c_{2}>0$ such that
$
f(C_{1}s)g(c_{1}s)\leq h(s)\leq f(C_{2}s)g(c_{2}s)
$
for the specified range of $s$.

The following theorem
is a special case of Theorems \ref{theorem:mainjump} and \ref{theorem:maindiff}.
Suppose that the fundamental solution
$p^0(t, x, y)$
of $\sL$ admits the following two-sided estimates:
\begin{equation}\label{e:1.7}
p^0(t, x, y)
\asymp t^{-d/\alpha} F(d(x, y)/t^{1/\alpha}),
\end{equation}
where either (i) $F(r) = \exp \left(-r^{\alpha/(\alpha-1)} \right)$ for
$\alpha\geq 2$;
or (ii) $F(r) = (1+r)^{-d-\alpha}$ with $\alpha>0$.
Note that for case (i), a priori $\alpha>1$ but in fact such estimates
can hold only if $\alpha\geq 2$ (see \cite[Page 1644]{bck}).
Case (i) typically corresponds to
diffusion processes. When
$$
\sL=\sum_{i,j=1}^d \frac{\partial}{\partial x_i} \left(a_{ij}(x)
\frac{\partial}{\partial x_j}\right)  \quad \hbox{with } \ \lambda^{-1}I_{d\times d} \leq (a_{ij}(x))\leq \lambda I_{d\times d}
\ \hbox{ on } \  \bR^d,
$$
where $I_{d\times d}$ denotes the $d\times d$ identity matrix,
it is known due to a result by Aronson \cite{Aro} that $\sL$ admits
 such estimates
with $\alpha=2$.
When $\sL$ is the Laplacian on a two-dimensional unbounded Sierpinski gasket,
it is shown by Barlow and Perkins in \cite{BP} that the
two-sided heat kernel estimates \eqref{e:1.7} of case (i) hold  $d=\log 3/\log 2$ and
$\alpha=d_w:=\log 5/\log 2$. Case (ii) typically corresponds to pure jump processes.
It is shown in \cite{CK1} that estimates \eqref{e:1.7} of case (ii) hold for
symmetric $\alpha$-stable-like process  on Alfhors $d$-regular space $E$ for $0<\alpha<2$.
See \cite{CKW} for examples with
$\alpha \ge 2$.

\smallskip

For $\beta \in (0, 1)$, define
\begin{align*}
H_{\le 1}(t,r)
=&\begin{cases}
 t^{\beta-1-\beta d/\alpha}, &  \text{if } d<2\alpha, \\
 t^{-1-\beta}\displaystyle \log \left(\frac{2t^{\beta }}{{r^\alpha}}\right),
 &  \text{if } d=2\alpha,\\
t^{-1-\beta}/r^{d-2\alpha}, &  \text{if } d>2\alpha,
\end{cases} \smallskip \\
H_{\ge 1}^{(j)} (t,r)
=& t^{2\beta-1}/r^{d+\alpha} ,   \smallskip \\
H_{\ge 1}^{(c)}(t,r)
=& t^{\beta-1-\beta d/\alpha}
\exp\Big(-(r^\alpha/t^{\beta})^{1/(\alpha-\beta)}\Big).
\end{align*}

\begin{thm}\label{T:1.4}
Suppose that $\sL$ is the generator of a Markov process and its
 fundamental solution $p^0(t, x, y)$
admits the two-sided estimates \eqref{e:1.7},
 and that $\{S_t, \P; t\geq 0\}$ is a $\beta$-stable subordinator with $0<\beta<1$.
 \begin{itemize}
\item[\rm (i)] Suppose $F(r)=(1+r)^{-d-\alpha}$ with $\alpha>0$.
 Then  the fundamental solution $q(t, x, y)$ for the time fractional Poisson equation \eqref{e:3.2}
satisfies
 \begin{align*}
q(t,x,y) \simeq
\begin{cases} H_{\le 1}(t,d(x,y))&\mbox{ if }  d(x,y)\leq  t^{\beta /\alpha},\\
H_{\ge 1}^{(j)}(t,d(x,y))
&\mbox{ if }  d(x,y) \geq t^{\beta /\alpha}.
\end{cases}
 \end{align*}
 \item[\rm (ii)] Suppose $F(r)=\exp(-r^{\alpha/(\alpha-1)})$ with
$\alpha\ge 2$.
Then  the fundamental solution $q(t, x, y)$ for the time fractional Poisson equation \eqref{e:3.2}
satisfies
 \begin{align*}
q(t,x,y)
\simeq H_{\le 1}(t,d(x,y))&\quad\mbox{ if }  d(x,y) \leq t^{\beta /\alpha} ,\\
q(t,x,y) \asymp
H_{\ge 1}^{(c)}
(t,d(x,y))&\quad\mbox{ if }  d(x,y) \geq t^{\beta /\alpha}.
 \end{align*}
\end{itemize}
\end{thm}

The rest of the paper is organized
as follows. Existence and uniqueness of
strong as well as weak
 solutions to the equation \eqref{e:3.2} are obtained in
Section \ref{S:2}.
Assuming further that $P_t^0$ has a density function $p^0(t,x,y)$ with respect to $\mu$, and that the subordinator $S$ is a special subordinator, we prove in Section \ref{section3} that
 the fundamental solution $q(t, x, y)$ for the time fractional Poisson equation and the fundamental solution $p(r, x, y)$ of the homogenous time fractional equation
 enjoy the following relation
\begin{align}
\label{e:tfder}
q(t, x, y) = \partial^{w^*}_t p (\cdot , x, y) (t)   \quad \hbox{ for  a.e.}\;\; t>0
\hbox{ and $\mu$-a.e. } x\not=y.
\end{align}
where  $w^*(x):=\nu^*(x, \infty)$ and $\nu^*$ is the L\'evy measure of the conjugate subordinator $S^*$ to $S$.
When  $\partial^{w }_t $ is the Caputo derivative of order $\beta \in (0, 1)$,
i.e., ${w }(s)=\frac{1}{\Gamma(1-\beta)} s^{-\beta}$, its corresponding subordinator
is the $\beta$-stable subordinator which
has $(1-\beta)$-stable subordinator as its conjugate.
Hence  $\partial^{w ^*}_t $ is the Caputo derivative of order $1-\beta$.
In this case,
that $\partial^{1-\beta}_t p(\cdot, x, y)(t)$
is the fundamental solution
for \eqref{e:3.2} has been established in \cite{EIK, EK, U, US} in some special cases of $\sL$,
and is called 
 Duhamel's principle in some literature.

In Section \ref{section4},
assuming that the distribution of subordinator is self-decomposable, we derive sharp estimates for the transition density function of the subordinator,
  which is of independent interest itself.
  Note that the self-decomposable distribution naturally occurs in scaling limits of random walks with  random waiting times.
See \cite[Corollary 3.8]{MSc}.
Using  estimates for the transition density function of the subordinator, we establish two-sided estimates for $q(t,x,y)$.
Full proofs of  the two-sided estimates for $q(t,x,y)$ is given in the appendix of this paper.

 Throughout the paper,
 for $a,b\in \bR$ we denote $a\wedge b:=\min\{a,b\}$ and $a\vee b:=\max\{a,b\}$.
For any measurable function $f$ we denote $f^+:=f\vee 0$.
 For a Banach space $(\bB, \| \cdot \|)$,
we also use $\| T \|$ to denote the operator norm for a linear operator $T:
(\bB, \| \cdot \|) \to (\bB, \| \cdot \|) $.

A constant  $c$ (without subscripts) denotes a strictly positive
constant,  whose value
is unimportant and which  may change from line to line.
Constants $c_0, c_1, c_2, \ldots$ with subscripts denote strictly positive
constants, and the labeling of the constants
starts anew in the statement of each result and the each step of its proof.

\section{Time fractional Poisson equations}\label{S:2}

Recall that throughout
 this paper,
$S=\{S_t, \P; t\geq 0\}$ is a driftless subordinator with infinite L\'evy measure $\nu$
and
starting from $0$; that is, $S_0=0$.
Recall that  $w(x):=\nu (x, \infty)$   and  that
$E_t =\inf\left\{ s>0: S_s >t\right\}$, $t>0$, is
the inverse subordinator of $S$.
 The assumption that the L\'evy measure $\nu$ is infinite (which is equivalent to $w(x)$
being unbounded) excludes  compound Poisson processes.
  Under this assumption, almost surely, $t\mapsto S_t$ is strictly increasing and hence $t\mapsto E_t$ is continuous.
Throughout this paper, we also
assume that \emph{the subordinator $S_r$
  has a bounded density function $\bar p(r, \cdot )$ for each $r>0$.}
This assumption holds, 
for example,
when Hartman and Wintner's condition is satisfied,
that is,
$$
\lim_{s \to \infty} \frac{\phi(s)}{\ln(1+ s)}=\lim_{s \to \infty}\frac{1} {\ln(1+ s)} \int_0^\infty (1-e^{-s x})\, \nu ( dx)  = \infty .
$$
See \cite[(74) on page 287]{HW}.

\smallskip

Since $\P (E_t \leq r) = \P(S_r >t)$, the following result
is from \cite{C}.

\begin{lemma} {\rm (\cite[Lemma 2.1]{C})} \label{T:2.1}
There is a Borel set $\sN \subset (0, \infty)$ having zero Lebesgue measure such that for every $t\in (0, \infty) \setminus \sN$, the inverse subordinator $E_t$ has a density function given by
\begin{equation}\label{e:2.1}
\frac{d}{dr} \P (E_t \leq r) =  \int_0^t w(t-s)\bar p (r, s) \,ds, \quad r>0.
\end{equation}
\end{lemma}

Let
\begin{equation}\label{e:2.2a}
G^{(S)}(t):=\int_0^\infty\bar p(r, t) \, dr
\end{equation}
 be  the  potential density of the subordinator $S$, and
denote by $G^{(S)}$ the 0-order resolvent of $S$; that is, for any  non-negative measurable function $f$ on $[0, \infty)$,
$$
G^{(S)}f(s)= \E \left[ \int_0^\infty f(s+S_t)\, dt\right]  = \int_0^\infty f(s+r) G^{(S)} (r)\, dr,\quad s\ge0.
$$
Since $S$ is transient, $
G^{(S)}{\bf 1}_{[0, T]} (0)<\infty$ for every $T>0$; see \cite[Theorem
35.4(v)]{Sat}.
In fact, by \cite[Proposition III.1]{Be},
 \begin{align}
\label{e:BeIII1}
G^{(S)} {\bf1}_{[0, T]} (0)\simeq 1/\phi (1/T) \quad  \text{ on }(0, \infty).
\end{align}

\begin{lemma}\label{L:1.20} Let
$$
w*G^{(S)}(t):=\int_0^t w(s)G^{(S)}(t-s)\,ds.
$$
Then $w*G^{(S)}(t) \le 1$ for all $t>0$, and $w*G^{(S)}(t) =1$  for a.e. $t>0$.
\end{lemma}

\pf
Let $\sN\subset (0, \infty)$ be the Borel set in Lemma \ref{T:2.1}.
Then, for every $t\in (0,\infty)\setminus \sN$,
\begin{align*} w*G^{(S)}(t)=&\int_0^t w(s)G^{(S)}(t-s)\,ds=\int_{0}^t w(t-s)G^{(S)}(s)\,ds\\
=& \int_0^\infty\int_0^t w(t-s) \bar p(r, s)\,ds\,dr = \bP (E_t <\infty) =1.\end{align*}

Furthermore, for each $t>0$, we take $t_n\notin \sN$ such that $t_n \downarrow t$ as $n\to \infty$.
 Then by Fatou's lemma and the right continuity of $w$, we conclude that
\begin{align*}
1= \lim_{n \to \infty}w*G^{(S)}(t_n) \ge \int_0^\infty  \liminf_{n \to \infty}
{\bf 1}_{\{ s\le t_n\}}w(t_n-s)G^{(S)}(s)\,ds=\int_0^t w(t-s)G^{(S)}(s)\,ds.
\end{align*} The proof is complete.
\qed

 The next proposition is a key step toward the existence and uniqueness of strong solutions for general time fractional Poisson equations.

\begin{prop}\label{P:2.3}
Suppose that $\{P^0_t; \, t\geq 0\}$ is a
 uniformly bounded and strongly continuous semigroup in some Banach space $(\bB, \| \cdot \|)$ and $T_0 \in (0,
 \infty)$.
 Then the following hold.
\begin{itemize}
\item[\rm (i)] For every $f(t, \cdot)$ with
  $\int_0^{T_0} \| f(t, \cdot )\| \,dt <\infty$,
\begin{equation}\label{e:2.2}
u(t, x) := \int_{s=0}^t \int_{r=0}^\infty P^0_r f (s, \cdot ) (x) \bar p (r, t-s)\, dr\, ds
\end{equation}
is well defined for a.e.\ $t\in [0,  T_0]$ as an element in $\bB$
such that $\int_0^{T_0} \| u(t, x)\| \, dt<\infty$.
Moreover,  for every $T\in (0,  T_0]$,
\begin{equation}\label{e:2.3}
\int_0^T w(T-t) u(t, x)\, dt = \int_0^T \E \left[ P^0_{E_{T-s}} f(s, \cdot ) (x)\right]
\, ds
\end{equation}
as elements in $\bB$
such that
\begin{equation}\label{e:2.3new}
\int_0^T w(T-t)  \| u(t, \cdot)\|\, dt \le M\int_0^T \| f(s, \cdot) \| \,ds,
\end{equation}
where  $M:=\sup_{t>0}\|P^0_t\| < \infty.$

\item[\rm (ii)] If there is some constant $K\in (0, \infty)$ so that $\| f(s, \cdot )\|\leq K$
for a.e.\ $s\in (0,  T_0]$, then $u(t, \cdot)$ is well defined
as an element in $\bB$ for every $t\in (0,  T_0]$ such that $\sup_{t\in (0,  T_0]} \| u(t, \cdot)\|
<\infty$, and
 $t\mapsto u(t,x)$ is continuous in $(\bB,\|\cdot\|)$ with $u(t, \cdot) \to 0$ in $\bB$ as $t\to 0$.

\item[\rm (iii)]
If $f(t, \cdot) \in \sD (\sL)$ for a.e. $t\in (0,  T_0]$ with
$\int_0^{T_0} \left( \|f(t, \cdot)\|+ \| \sL f (t, \cdot) \| \right)\, dt <\infty$,
then $u(t, \cdot) \in \sD(\sL)$ for a.e. $t\in (0,  T_0]$ with $\int_0^{T_0}
\| \sL u(t, \cdot )\| \,dt <\infty$.
Moreover, if there is some constant $K\in (0, \infty)$ so that
$ \| f(s, \cdot )\|+ \| \sL f(s, \cdot )\|\leq K$
for a.e.\ $s\in (0,  T_0]$, then
 $t\mapsto \sL u(t,x)$ is continuous in $(\bB,\|\cdot\|)$.
\end{itemize}
\end{prop}

   \pf   Throughout the proof,   let $M=\sup_{t>0}\|P^0_t\|<\infty$ and $T_0 \in (0,  \infty)$.

\noindent (i)
For given $f(t, \cdot)$ with
  $\int_0^{T_0} \| f(t, \cdot )\| \,dt <\infty$, by the Fubini theorem,
\begin{equation}\label{e:p2.31}\begin{split}
 \int_{t=0}^{T_0} \int_{s=0}^t \int_{r=0}^\infty
\| P^0_r f(s, \cdot) \| \, \bar p(r, t-s) \,dr\, ds \,dt
&\leq M \int_{s=0}^{T_0} \int_{t=s}^{T_0} \int_{r=0}^\infty
 \| f(s, \cdot )\| \, \bar p(r, t-s) \,dr \,
 dt \,ds\\
&= M \int_{s=0}^{T_0} \int_{r=0}^\infty
 \| f(s, \cdot )\| \, \P(S_r \leq {T_0}-s)\, dr\, ds \\
&\leq M \int_0^{T_0} \| f(s, \cdot) \|\,ds \, G^{(S)}{\bf1}_{[0, {T_0}]} (0).
\end{split}\end{equation}
This establishes the well-definedness of $u(t, x)$
for a.e.\ $t\in (0,  T_0]$ with
$$
\int_0^{T_0} \| u(t, \cdot ) \|\, dt
\leq M \int_0^{T_0} \| f(s, \cdot) \|\, ds \, G^{(S)}{\bf 1}_{[0, {T_0}]} (0) <\infty.
$$

We now show that both sides of \eqref{e:2.3} are well defined as an element in $\bB$.
Note that, for every $T\in (0,  T_0]$,
$$
 \int_0^T \E \| P^0_{E_{T-s}} f (s, \cdot) \|\, ds \leq M \int_0^T
\| f(s, \cdot) \| \,ds <\infty.
$$
 This shows that $ \int_0^T  \E [P^0_{E_{T-s}} f (s, \cdot)]  \,ds$
is absolutely integrable and thus it is well defined as an element in $\bB$.
Similarly,
by the definition of $u(t, x)$ and the Fubini theorem,
\begin{align*}
\int_0^T w(T-t) \| u(t, \cdot ) \| \,dt
&\leq  \int_{t=0}^T w(T-t)  \left(\int_{s=0}^t \int_{r=0}^\infty
\| P^0_r f(s, \cdot )\| \, \bar p(r, t-s) \,dr \,ds \right)\, dt \\
&\leq M\int_{s=0}^T  \left(\int_{t=s}^T\int_{r=0}^\infty
w(T-t) \bar p(r, t-s) \,dr\,dt \right) \| f(s, \cdot ) \|\, ds \\
&= M\int_{s=0}^T \left( \int_{t=s}^T
w(T-t) G^{(S)}(t-s)\, dt  \right) \| f(s, \cdot ) \|\, ds \\
&=M \int_0^T \| f(s, \cdot) \| \,ds<\infty,
\end{align*}
where in the last equality we used Lemma \ref{L:1.20}.
This shows that $ \int_0^T w(T-t) u(t, \cdot )\, dt$
is absolutely integrable and thus
 it is well defined as an element in $\bB$.

Next, we verify that the identity \eqref{e:2.3} holds.
By  \eqref{e:2.1}, Lemma \ref{L:1.20} and the Fubini theorem,
\begin{align*} \int_0^T \E \left[ P^0_{E_{T-s}} f(s, \cdot ) (x ) \right]\, ds &= \int_{s=0}^T \int_{r=0}^\infty P^0_r f(s, \cdot) (x)\, d_r \P (E_{T-s} \leq r)  \,ds \\
&= \int_{s=0}^T \int_{r=0}^\infty P^0_r f(s, \cdot) (x) \left( \int_{t=0}^{T-s}
w(T-s-t) \bar p(r, t)\, dt \right)\,dr \,ds \\
&=   \int_{s=0}^T \int_{r=0}^\infty P^0_r f(s, \cdot) (x)  \left( \int_{t=s}^{T}
w(T-t) \bar p(r, t-s) \,dt \right)   \,dr \, ds \\
&=  \int_{t=0}^{T} \left( \int_{s=0}^t   \int_{r=0}^\infty P^0_r f(s, \cdot) (x)
 \bar p(r, t-s) \,dr \,ds \right)  w(T-t)    \,dt \\
&= \int_{0}^{T}  w(T-t) u(t, x)\, dt.
\end{align*}
This establishes \eqref{e:2.3}.

\noindent
(ii)
If $\| f(s, \cdot )\|\leq K<\infty$
for a.e.\ $s\in (0,  T_0]$, then, by
 \eqref{e:BeIII1}
we have
that for every $t\in (0,  T_0]$,
\begin{equation}\label{e:2.4}\begin{split}
\int_{s=0}^t \int_{r=0}^\infty \| P^0_r f(s, \cdot) \| \, \bar p(r, t-s) \,dr\, ds
&\leq  M\, \int_{s=0}^t \int_{r=0}^\infty \|  f(s, \cdot) \| \, \bar p(r, t-s)\, dr\, ds
  \\
&\leq M \int_0^t \|f(s, \cdot) \| \, G^{(S)}(t-s)\, ds  \\
&\leq
M K\, G^{(S)} {\bf1}_{[0 , t]}(0)
\le M K\,  \frac{c_0}{\phi(1/t)}.
\end{split}\end{equation}
So
$u(t, x):=\int_{s=0}^t \int_{r=0}^\infty P^0_r f(s, \cdot)\bar p(r, t-s)\,dr \, ds$
is well defined as an element in $\bB$ and that for
every $t\in (0,  T_0]$, $\sup_{t\in (0,  T_0]} \| u(t, \cdot)\| \le MK {c_0}/{\phi(1/T_0)}
<\infty$.
Since $\lim_{r\uparrow \infty} \phi (r)=\nu(0,\infty)=\infty$,
we have $\lim_{t\downarrow 0} \| u(t, \cdot )\| = 0$.

\medskip

We now show the continuity of $t\mapsto u(t,x)$ in $(\bB,\|\cdot\|)$.
Indeed, for any $0<t_1\le t_2$,
\begin{equation}\label{e:cont1}
\begin{split}
u(t_2,x)-u(t_1,x)
=&\int_{s=t_1}^{t_2} \int_{r=0}^\infty P^0_r f (s, \cdot ) (x) \bar p (r, t_2-s)\, dr\, ds\\
&+\int_{s=0}^{t_1} \int_{r=0}^\infty P^0_r f (s, \cdot ) (x) (\bar p(r,t_2-s)-\bar p (r, t_1-s))\, dr\, ds. \end{split}\end{equation}
Thus,
\begin{equation}\label{e:cont2}\begin{split}
\|u(t_2, \cdot)-u(t_1, \cdot)\|
\le &\int_{s=t_1}^{t_2} \int_{r=0}^\infty \| P^0_r f (s, \cdot )\|\ \bar p (r, t_2-s)\, dr\, ds\\
&+\int_{s=0}^{t_1} \int_{r=0}^\infty  \|P^0_r f (s, \cdot ) \| |\bar p(r,t_2-s)-\bar p (r, t_1-s)|\, dr\, ds\\
\le &M K\int_{s=t_1}^{t_2} \int_{r=0}^\infty \bar p (r, t_2-s)\, dr\, ds\\
&+MK\int_{s=0}^{t_1} \int_{r=0}^\infty  |\bar p(r,t_2-s)-\bar p (r, t_1-s)|\, dr\, ds\\
= &M K
 G^{(S)} {\bf1}_{[0, t_2-t_1]}(0)
+MK\int_{s=0}^{t_1} \int_{r=0}^\infty  |\bar p(r,t_2-t_1+s)-\bar p (r, s)|\, dr\, ds.\end{split}\end{equation}
 Using \eqref{e:BeIII1}
and the fact $\lim_{r\uparrow \infty} \phi (r)=\nu(0,\infty)=\infty$ in the first term, and the
$L^1$-continuity of the translation of the integrable function
$(r, s) \mapsto \bar p (r, s)$
on $(0, \infty) \times (0,  T_0]$ on the second term, we see that both terms  go to zero as $t_2-t_1 \to 0$.
 Thus  the desired assertion follows.

\noindent
(iii) Suppose that $f(t, \cdot) \in \sD (\sL)$ for a.e. $t\in (0,  T_0]$ with
$\int_0^{T_0} \left( \|f(t, \cdot)\|+ \| \sL f (t, \cdot) \| \right)\, dt <\infty$.
Note that for a.e.\ $s\in [0,  T_0]$ and every $r>0$, $\sL P^0_r f(s,
\cdot) =P^0_r \sL f (s, \cdot)$ and $\| \sL P^0_r f(s, \cdot)\|  =
\| P^0_r \sL f (s, \cdot) \|\leq M \| \sL f(s, \cdot)\|.$
By the similar calculation as that in  \eqref{e:p2.31},
we have
$$
\int_{t=0}^{T_0}\int_{s=0}^{T_0} \int_{r=0}^\infty \| \sL P^0_r f(s,
\cdot ) \| \, \bar p(r, t-s)\, dr\, ds \,dt \leq M \int_0^{T_0} \|
\sL f(s, \cdot ) \|\,ds \, G^{(S)}{\bf1}_{[0,  T_0]}(0) <\infty.
$$
Thus, by \eqref{e:2.2}
and the closed graph theorem,
$u(t, \cdot ) \in \sD (\sL)$ for a.e.\
$t\in [0,  T_0]$ such that
\begin{align}\label{e:intlbd}
\int_0^{T_0}  \| \sL u(s, \cdot )\|\, ds \leq c
 \int_0^{T_0}  \| \sL f(s, \cdot ) \|\,ds \,
G^{(S)}{\bf1}_{[0,  T_0]}(0)
 <\infty.
\end{align}

Next, we suppose that $\| f(s, \cdot )\|+\|  \sL f(s, \cdot )\|\leq K<\infty$
for a.e.\ $s\in (0,  T_0]$. Then,
by the closed graph theorem, the Riemann sum approximation
and the fact that $\sL P^0_r f (s, \cdot )=P^0_r \sL f (s, \cdot )$ for a.e.\ $s\in (0,  T_0]$ and any $r>0$,
we have that for all $0\leq T \le T_0$.
\begin{align}\label{e:cgR}
\int_0^T \sL u(t, \cdot) (x)\, dt
=& \int_{t=0}^T \int_{s=0}^t
\int_{r=0}^\infty  P^0_r \sL f(s, \cdot) (x)
\bar p (r, t-s)\, dr \,ds \,dt.
\end{align}
Using \eqref{e:cgR}, we follow exactly the
same arguments  these in \eqref{e:cont1} and \eqref{e:cont2} and get that for any $0<t_1\le t_2$,
\begin{align*}
\| \sL u(t_2,\cdot)- \sL u(t_1,\cdot)\|
\le &\int_{s=t_1}^{t_2} \int_{r=0}^\infty \| P^0_r  \sL f (s, \cdot )\|\ \bar p (r, t_2-s)\, dr\, ds\nn\\
&+\int_{s=0}^{t_1} \int_{r=0}^\infty  \|P^0_r  \sL f (s, \cdot ) \| |\bar p(r,t_2-s)-\bar p (r, t_1-s)|\, dr\, ds\nn\\
\le &M K
G^{(S)} {\bf1}_{[0, t_2-t_1]}(0)
+MK\int_{s=0}^{t_1} \int_{r=0}^\infty  |\bar p(r,t_2-t_1+s)-\bar p (r, s)|\, dr\, ds.\end{align*}
Therefore,  $t\mapsto \sL u(t,x)$ is continuous in $(\bB,\|\cdot\|)$.
This completes the proof of the proposition. \qed

\begin{remark} \rm If we assume
$$
 \int_0^{T_0} \| f(s, \cdot)\| G^{(S)} (t-s)\, ds<\infty  \,
\hbox{ with } \,
\lim_{t\to 0} \int_0^t \| f(s, \cdot)\| G^{(S)} (t-s) \,ds=0
$$
in place of $\|f(t, \cdot)\| \leq K$ for a.e. $t\in (0,  T_0]$ in Proposition \ref{P:2.3}(ii),
then, by the same proof, $u(t, \cdot)$ is well defined
as an element in $\bB$ for every $t\in (0,  T_0]$ with $u(t, \cdot) \to 0$ in $\bB$ as $t\to 0$.
\end{remark}

\bigskip

Recall from  \cite[(2.3)]{C} that
\begin{equation}\label{e:3.6}
\int_0^\infty e^{-\lambda x} w(x) \,dx = \frac{\phi (\lambda)}{\lambda}
 \quad \hbox{for } \lambda >0.
 \end{equation}

\bigskip

\noindent {\bf Proof of Theorem \ref{T:mains}.}\,\, (a) {\bf(Existence)}\,\,
The case that $f\equiv 0$ is proved in \cite[Theorem 2.3]{C}. Thus
 by the linearity  we only need to show the existence of the solution when
   $g \equiv 0$, which we will assume throughout the proof.

Let $u(t, x)$ be defined by \eqref{e:2.2}.
Part (i) has been proved in Proposition
\ref{P:2.3}(ii)
as $\| f (t, \cdot )\|\leq K$ for a.e. $t\in (0,  T_0]$,
while part (ii) has been established in Proposition \ref{P:2.3}(iii).

By \eqref{e:intlbd}, using the closed graph theorem and the Riemann sum approximation,
we have that for every $T\in (0,  T_0]$,
$$
\int_0^T \sL u(t, \cdot) (x) \,dt
= \int_{t=0}^T \int_{s=0}^t
\int_{r=0}^\infty \sL P^0_r f(s, \cdot) (x)
\bar p (r, t-s)\, dr \,ds \,dt.
$$
Thus, for every $T\in (0,  T_0]$, using the Fubini theorem and
the integration by parts,
\begin{align*}
 \int_0^T \sL u(t, \cdot) (x) \,dt
&=\int_{s=0}^T \int_{r=0}^\infty \sL P^0_r f(s, \cdot) (x)
\left( \int_{t=s}^T  \bar p (r, t-s) \,dt \right) \,dr \,ds  \\
&=\int_{s=0}^T  \left( \int_{r=0}^\infty \frac{d}{dr} P^0_r f(s,
\cdot) (x)
\P(S_r \leq T-s) \,dr \right)\, ds  \\
&=\int_{s=0}^T  \left(  P^0_r f(s, \cdot) (x) \P(S_r \leq T-s)
\Big|_{r=0}^\infty
 -\int_{r=0}^\infty   P^0_r f(s, \cdot) (x)\,
d_r \P(S_r \leq T-s) \right) ds  \\
&= \int_{s=0}^T  \left( -f(s, x) +\int_{r=0}^\infty   P^0_r f(s,
\cdot) (x)\, d_r \P(E_{T-s} \leq r ) \right) \,ds .
\end{align*}
Hence,
\begin{equation}\label{e:2.11}
 \int_0^T \left( f(t, x) + \sL u (t, \cdot) (x) \right) dt
= \int_0^T \E \left[ P^0_{E_{T-s}} f (s,\cdot ) (x) \right] ds =
\int_0^T w(T-t) u(t, x)\, dt,
\end{equation}
where the last equality is due to \eqref{e:2.3}. The proof of (iii)
is finished.

\medskip

(b) {\bf (Uniqueness)}
Suppose that $u$ and $v$ are two (strong)
solutions of
$\partial^w_t u (t,x)= \sL u(t,x) + f(t, x)$
 on $(0,  T_0]\times E$ with $u(0, x)=g(x)$.
Then $h:=u-v$ is
a (strong) solution of
\begin{equation}\label{e:2.C21}
 \partial^w_t h(t,x)  =
\sL h (t,x)
\quad \hbox{on } (0,  T_0]\times E \
\hbox{ with } h (0, \cdot )=0.
\end{equation}
By the definition of
 the strong solution, we know that $t \mapsto h(t, \cdot)$ is continuous on $[0,  T_0]$
with $\lim_{t\to 0}h(t, \cdot)=0$ and $K_1:=\sup_{t\in [0,  T_0]} \| h(t, \cdot )\|<\infty$;
moreover,
\begin{equation}\label{e:2.12}
h(t, \cdot) \in \sD (\sL) \hbox{ for a.e. } t\in (0,  T_0] \ \hbox{ with } \  \int_0^{T_0} \| \sL h (t, \cdot)\| \,dt <\infty.
\end{equation}
Define
\begin{equation}\label{e:dfbarf}
\bar f (t, x):=\int_{(t,T_0+t]} h(T_0+t -s, x)\,
\nu (ds), \quad t \ge 0.
\end{equation}
Clearly,
 $\| \bar f (t, \cdot )\| \leq K_1 w(t)$ for $t>0$,   and so by \eqref{e:3.6} we have that
 for every $\lambda >0$,
\begin{equation}\label{e:2:C1}
\int_0^\infty e^{-\lambda t} \|\bar f (t, \cdot )\| \,dt \leq K_1 \int_0^\infty e^{-\lambda t} w(t)\, dt
=K_1\phi (\lambda)/\lambda.
\end{equation}
Moreover, since
$\int_{0}^\infty e^{-\lambda t} \,\nu(dt) = \phi'(\lambda)\lambda \le \phi (\lambda)$,  by \eqref{e:2.12}
and the Fubini theorem,
we have that for every $\lambda >0$
\begin{equation}\label{eq:biefiqlc}
\begin{split}
 \int_{0}^\infty e^{-\lambda t} \int_{ t}^{T_0+t}\|
\sL  h(T_0+t -s, \cdot) \|\,\nu (ds) \, dt
&=
\int_{0}^\infty \int_{(s-T_0)^+ }^{s} e^{-\lambda t} \|   \sL  h(T_0+t -s, \cdot) \|\,dt
 \,\nu (ds)   \\
&=
\int_{0}^\infty  \int_{ (T_0-s)^+}^{T_0} e^{-\lambda (r+s-T_0)} \|   \sL  h(r, \cdot) \| \,dr \,
\nu (ds)  \\
&\le
e^{\lambda T_0} \left(\int_{0}^\infty e^{-\lambda s}\,
\nu (ds) \right) \left(\int_{0}^{T_0}  \|   \sL  h(r, \cdot) \|\,dr \right)   \\
&\le   e^{\lambda T_0} \phi(\lambda) \int_0^{T_0} \|   \sL  h(t, \cdot) \|\,
dt<\infty.
\end{split}
\end{equation}
Hence, by the closed graph theorem and the Riemann sum approximation, we conclude that
$\bar f (t, \cdot) \in \sD(\sL)$
for a.e. $t>0$, and
$$\sL \bar f (t, \cdot ) (x)
= \int_{(t,T_0+t]}\sL h(T_0+t -s, \cdot ) (x)\,
\nu (ds)  \quad \hbox{for a.e. } t > 0
$$
 with
\begin{equation}\label{e:2.16}
\int_0^\infty e^{-\lambda t}\|  \sL \bar f(t, \cdot) \| \,dt \leq e^{\lambda T_0} \phi(\lambda) \int_0^{T_0} \|   \sL  h(t, \cdot) \|\,dt<\infty.
\end{equation}

Let
\begin{equation}\label{e:dfbarh}
\bar h (t, x) := \int_{s=0}^t \int_{r=0}^\infty P^0_r \bar f (s, \cdot ) (x) \bar p (r, t-s)\, dr\, ds,   \quad t>0,
\end{equation}
and $M:=\sup_{t>0}\|P^0_t\|
<\infty.$
 By \eqref{e:2:C1}, \eqref{e:2.16} and Proposition \ref{P:2.3},  for a.e.\ $t\in [0,\infty)$, $\bar h(t, \cdot)$ is well defined as an element in $\sD(\sL) \subset \bB$ with
$\int_0^T \left( \| \sL \bar h(t, \cdot )\|+ \| \bar h(t, \cdot) \|  \right)dt <\infty$
for every $T>0$.
Moreover, by \eqref{e:2:C1} and the Fubini theorem,
for every $\lambda >0$
\begin{align}
\int_0^\infty e^{-\lambda t}  \| \bar h (t, \cdot)\| \, dt
&\le
\int_{t=0}^{\infty}  e^{-\lambda t} \int_{s=0}^t \int_{r=0}^\infty
\| P^0_r \bar f(s, \cdot) \| \, \bar p(r, t-s) \,dr\, ds \,dt \nn\\
&\leq  M \int_{s=0}^{\infty} \int_{t=s}^{\infty}  e^{-\lambda t}\int_{r=0}^\infty
 \| \bar f(s, \cdot )\| \, \bar p(r, t-s) \,dr \,
 dt \,ds  \nn\\
 &=  M  \int_{s=0}^{\infty}  e^{-\lambda s} \|\bar f(s, \cdot )\| \, \left(\int_{t=s}^{\infty}  e^{-\lambda (t-s)}\int_{r=0}^\infty
 \bar p(r, t-s) \,dr \,
 dt \right) \,ds \nn\\
  &=  \frac{M}{\phi(\lambda)} \int_{s=0}^{\infty}  e^{-\lambda s} \|\bar f(s, \cdot )\| \,ds
	 \le  MK_1 /{\lambda}   <\infty.\label{e:lambdabd}
\end{align}
Here in the  last equality we used the fact that
\begin{equation}\label{e:2.20}
\int_{t=s}^{\infty}  e^{-\lambda (t-s)}\int_{r=0}^\infty
 \bar p(r, t-s) \,dr \, dt
= \int_{r=0}^\infty \bE e^{-\lambda S_r}\, dr = \int_0^\infty e^{-r \phi (\lambda )}\, dr = 1/\phi (\lambda),
\end{equation}
and the last inequality follows from \eqref{e:2:C1}.
In view of all the estimates above,
we have by the same argument as that for \eqref{e:2.11} (where the boundedness of $t\mapsto \|f (t, \cdot )\|$ is not needed) that for every $T>0$,
\begin{equation}\label{e:2.21}
\int_0^T \left( \bar f(t, x) + \sL \bar h (t, \cdot) (x) \right) \,dt  =
\int_0^T w(T-t) \bar h (t, x)\, dt.
\end{equation}

Now, we extend $h(t,x)$ on $(0,\infty)\times E$ by defining $h(t, x)= \bar h (t-T_0, x)$ for $t>T_0$.
Then, since for all $T>T_0$
\begin{equation}\label{e:TT0}\begin{split}
\int_0^{T-T_0}
\bar f (t, x)\,dt&=
\int_0^{T-T_0} \int_{(t,T_0+t]} h(T_0+t -s, x)
\,\nu (ds)\,dt \\
&=
\int_{(0,T]} \int_{0 \vee (s-T_0)}^{s \wedge (T-T_0)}
h(T_0+t -s, x) \,dt \,\nu (ds)\\
& =
\int_{(0,T]} \int_{0 \vee (T_0-s)}^{T_0 \wedge (T-s)}
h(t , x) \,dt \,\nu (ds)\\
&=\int_{0}^{T_0}\left(\int_{(T_0-t, T-t]} \nu(ds)\right)h(t,x)\,dt\\
&=
\int_{0}^{T_0}\left(w(T_0-t)-w(T-t)\right)h(t,x)\,dt,
\end{split}\end{equation}
 we have that, by \eqref{e:2.C21}, \eqref{e:TT0} and  \eqref{e:2.21}, for all $T>T_0$
\begin{align*}
 \int_0^T  \sL   h (t, \cdot)(x)\,dt
&=\int_0^{T_0} w(T_0-t)h(t,x)\,dt + \int_{T_0}^T  \sL \bar h (t-T_0, x)\,dt\\
&=
\int_0^{T_0} w(T-t)h(t,x)\,dt+\int_{0}^{T_0}\left(w(T_0-t)-w(T-t)\right)h(t,x)\,dt\\
&\quad
 +  \int_{0}^{T-T_0}  \sL \bar h (t, x)\,dt\\
&= \int_0^{T_0} w(T-t)h(t,x)\,dt+\int_0^{T-T_0}
(\bar f (t, x)dt+\sL \bar h (t, x))\,dt\\
&= \int_0^{T_0} w(T-t)h(t,x)\,dt+\int_0^{T-T_0}
 w(T-T_0-t) \bar h(t, x) \,dt\\
 &= \int_0^{T_0} w(T-t)h(t,x)\,dt+\int_{T_0}^{T}w(T-t)
 \bar h (t-T_0, x)\,dt \\
&=\int_0^{T} w(T-t)h(t,x)\,dt.
\end{align*}
This shows that
$h$ solves
\begin{equation}\label{e:2.23}
 \partial^w_t   h(t, x) = \sL   h (t, x) \quad \hbox{for all $t>0$ with } h(0, x)=0.
\end{equation}
To show that $h(t, \cdot)=0$ to a.e. $t>0$, first note that
for every $\lambda >0$,
\begin{equation}
\label{e:lambdabd1}\begin{split}
&\int_0^\infty e^{-\lambda t}  \| \sL h (t, \cdot)\| \, dt  \\
&=\int_0^{T_0}  e^{-\lambda t}   \| \sL h (t, \cdot)\| \, dt +
\int_{T_0}^\infty e^{-\lambda t}  \| \sL \bar h (t-T_0, \cdot)\| \, dt   \\
& \le \int_0^{T_0} \| \sL h (t, \cdot)\| \, dt +
\int_{0}^\infty e^{-\lambda t}  \| \sL \bar h (t, \cdot)\| \, dt  \\
& \le \int_0^{T_0} \| \sL h (t, \cdot)\| \, dt + \int_{t=0}^{\infty}  e^{-\lambda t} \int_{s=0}^t \int_{r=0}^\infty
\| \sL P^0_r \bar f(s, \cdot) \| \, \bar p(r, t-s) \,dr\, ds \,dt \\
&\leq \int_0^{T_0} \| \sL h (t, \cdot)\| \, dt + M \int_{s=0}^{\infty} \int_{t=s}^{\infty}  e^{-\lambda t}\int_{r=0}^\infty
 \| \sL \bar f(s, \cdot )\| \, \bar p(r, t-s) \,dr \,
 dt \,ds \\
 &=\int_0^{T_0} \| \sL h (t, \cdot)\| \, dt + M \int_{s=0}^{\infty}  e^{-\lambda s} \|\sL \bar f(s, \cdot )\| \, \left(\int_{t=s}^{\infty}  e^{-\lambda (t-s)}\int_{r=0}^\infty
 \bar p(r, t-s) \,dr \,
 dt \right) \,ds \\
  &= \int_0^{T_0} \| \sL h (t, \cdot)\| \, dt + \frac{M}{\phi(\lambda)} \int_{s=0}^{\infty}  e^{-\lambda s} \|\sL \bar f(s, \cdot )\| \,ds \\
&\leq (1+Me^{\lambda T_0}) \int_0^{T_0} \| \sL h(t, \cdot) \| \,dt <\infty,
\end{split}\end{equation}
where \eqref{e:2.20} was used the last equality,
and  \eqref{e:2.12} and \eqref{e:2.16}
were used in the last  inequality.
Now repeating the proof of the uniqueness part in \cite[Theorem 2.3]{C}, we can show that $h(t, \cdot )=0$ for a.e. $t>0$. Indeed,
let $H(\lambda, x):=\int_0^\infty e^{-\lambda t} h(t, x)\, dt$, $\lambda >0$,
 be the Laplace transform of $t\mapsto h(t, x)$. By \eqref{e:lambdabd} for every $\lambda >0$, $H(\lambda, \cdot)\in \bB$
with $\| H(\lambda, \cdot)\| <\infty$. By the closed graph theorem, the Riemann sum approximation  and \eqref{e:lambdabd1},
for each $\lambda >0$, $H(\lambda, \cdot )\in \sD (\sL)$ with
$$
\sL H(\lambda, \cdot) = \int_0^\infty e^{-\lambda t} \sL h(t, \cdot)\, dt
\quad \hbox{and} \quad
\|\sL H(\lambda, \cdot) \|\leq \int_0^\infty e^{-\lambda t}  \| \sL h(t, \cdot) \| \,dt < \infty.
$$
Taking the Laplace transform
 in $t$ on both sides of \eqref{e:2.23} yields
$$H(\lambda, x)  \frac{\phi(\lambda)}{\lambda}=
H(\lambda, x)
 \int_0^\infty e^{-\lambda s} w(s) \,ds
= \frac{1}{\lambda} \int_0^\infty  e^{-\lambda t} \sL h(t, x) \,dt = \frac{\sL H(\lambda, x)}{\lambda}.
$$
Thus
$
(\phi (\lambda)-\sL ) H(\lambda, x) = 0$ for every  $\lambda >0.
$
Since $\sL$ is the infinitesimal generator of a
 uniformly bounded and strongly continuous
semigroup $\{P^0_t, t\geq 0\}$ in Banach space $\bB$,
for every $\alpha >0$, the resolvent $G^0_\alpha = \int_0^\infty e^{-\alpha t} P^0_t\,dt$ is well defined and is the inverse to $\alpha-\sL$.
Hence $H(\lambda, \cdot)=0$ in $\bB$ for every $\lambda >0$. By the uniqueness of the Laplace transform, we conclude that $h(t, \cdot )=0$ for a.e. $t>0$.

Since $t\mapsto h(t, \cdot)$ is continuous in $\bB$ for $t\in [0,  T_0]$, we conclude that $h(t, \cdot)=0$
in $\bB$ for every $t\in [0,  T_0]$. This establishes the uniqueness.
\qed

\medskip

\begin{remark}\label{R:2.4} \rm
If we assume that there is some positive constant $K$ so that $ \|
\sL f(t, \cdot )\|\leq K$ for a.e.\ $t\in (0,  T_0]$ instead of
$\int_0^{T_0}  \| \sL f(t, \cdot) \|\,  dt<\infty$ in Theorem
\ref{T:mains}, then its conclusion (ii) can be strengthened to
\begin{description}
\item{\rm (ii')} $u(t, \cdot)\in \sD (\sL)$ for every $t\in (0,  T_0]$
and $\sup_{t\in [0,  T_0]} \| \sL u(t, \cdot ) \|<\infty$.
\end{description}
Its proof is similar to the last part of that for Proposition
\ref{P:2.3}.
\end{remark}

\medskip

To study weak solutions for general time fractional  Poisson
equations, we will assume
that $(\sL, \sD (\sL))$ is the
infinitesimal generator of a  (symmetric) Dirichlet
 form on $L^2(E; \mu)$. Its
associated transition semigroup is denoted as $\{P^0_t; t\geq 0\}$,
which is a strongly continuous contraction semigroup in $L^2(E;
\mu)$. Denote by $\langle \cdot,\cdot\rangle$ the inner product of
$L^2(E;\mu)$.

\begin{prop} \label{T:2.5}
Suppose that $f(t, x)$ is a function on $(0,  T_0]\times E$ with
$ \int_0^{T_0} \| f(t, \cdot )\|_{L^2(E; \mu)}\, dt <\infty$
 for some $T_0\in (0,
 \infty)$.
Then the following hold.

 \noindent
{\rm (1)}
$u(t,x):=\int_{s=0}^t \int_{r=0}^\infty P^0_r f(s,
\cdot)\bar p(r, t-s)\,dr \, ds$
is well defined as an element in $L^2(E; \mu)$ for a.e. $t\in (0,  T_0]$,
and it is a  weak solution of
$\partial^w_t u = \sL u + f(t, x)$ on $(0,  T_0]\times E$ with $u(0,
x)=0$ in the following sense:

\begin{description}
\item{\rm (i)}
$u(t, \cdot)\in L^2(E; \mu)$ for a.e. $t\in (0,  T_0]$ with
$\int_0^{T_0} \| u(t, \cdot )\|_{L^2(E; \mu)}\, dt <\infty$.
Moreover, for every $T\in (0,  T_0]$,
$$
\int_0^T w(T-t)  \| u(t, \cdot)\|_{L^2(E; \mu)}\, dt <\infty
$$
and so $\int_0^T w(T-t) u(t, \cdot)\, dt$ is well defined as an element in $L^2(E; m)$.

\item{\rm (ii)} For every $\varphi \in \sD (\sL)$ and $T\in (0,  T_0]$,
\begin{equation}\label{e:2.7}
\Big\< \int_0^T w(T-s) u(s, \cdot)\,ds, \varphi \Big\>
=\int_0^T\<f(s, \cdot), \varphi \> \,ds + \int_0^T \< u(s, \cdot), \sL
\varphi \> \,ds
\end{equation}
\end{description}

\noindent
{\rm (2)}  If we further assume that
\begin{equation}\label{e:2.7n}
\| f(t, \cdot )\|_{L^2(E; \mu)} \leq K<\infty \ \hbox{ for a.e. } t\in (0,  T_0],
\end{equation}
then the conclusion {\rm (i)} of {\rm (1)}
can be strengthened to
\begin{description}
\item{\rm (i')} $u(t, \cdot)$ is well defined as an element in $L^2(E; \mu)$ for every $t\in (0,  T_0]$
satisfying that $t\mapsto u(t,x)$ is continuous in
$L^2(E; \mu)$,
$$\sup_{t\in (0,  T_0]} \| u(t, \cdot )\|_{L^2(E; \mu)}  <\infty
\quad \hbox{and} \quad u(t, \cdot) \to 0 \hbox{ in } L^2(E; \mu) \hbox{  as  } t\to 0.
$$
If in addition, $f(t, \cdot) \in \sD (\sL)$ for a.e. $t\in (0,  T_0]$ with
$\| \sL f(t, \cdot)\|_{L^2(E; \mu)} \leq K <\infty$ for a.e. $t\in (0,  T_0]$,
then $t\mapsto \sL u(t,x)$ is  continuous in $L^2(E; \mu)$.
\end{description}
\end{prop}

\pf
 \noindent
{\rm (1)}
  Let $u(t, x):= \int_{s=0}^t
\int_{r=0}^\infty P^0_r f(s, \cdot)\bar p(r, t-s)\,dr \, ds $. By
Proposition \ref{P:2.3}(i), $u(t, x)$ is well defined for a.e.\ $t\in [0,  T_0]$ as an element in
$L^2(E; \mu)$  having $\int_0^{T_0} \| u(t, x)\|_{L^2(E; \mu)} \, dt<\infty$.
It follows from \eqref{e:2.3new} that for every $T\in (0,  T_0]$,
\begin{align}
\label{e:Twf}
\int_0^T w(T-t)  \| u(t, \cdot)\|_{L^2(E; \mu)}\, dt \le M
\int_0^T \| f(s, \cdot) \|_{L^2(E; \mu)} \,ds<\infty
\end{align}
 with $M=\sup_{t>0}\|P^0_t\| $.
Hence $\int_0^T w(T-t)   u(t, \cdot)\, dt $ is well defined as element in $L^2(E;\mu)$.

\medskip

For every $\varphi \in \sD (\sL)$ and $T\in (0,  T_0]$,  by \eqref{e:2.3}, \eqref{e:Twf}, the Fubini theorem and the symmetry of $P^0_r$,
\begin{equation}\label{e:2.8}\begin{split}
 \left\< \int_0^T w(T-t) u(t, \cdot)\,dt, \varphi \right\>
&= \left\< \int_0^T \E P^0_{E_{T-s}} f(s, \cdot)\, ds, \varphi \right\> \\
&= \int_0^T \< f(s, \cdot), \E P^0_{E_{T-s}} \varphi \> \,ds\\
&= \int_0^T \< f(s, \cdot), \varphi\> \,ds + \int_0^T \left\< f(s,
\cdot),
\E \int_0^{E_{T-s}} P^0_t \sL \varphi \, dt \right\> \,ds \\
&= \int_0^T \< f(s, \cdot), \varphi\>\, ds + \int_0^T \left\< \E
\int_0^{E_{T-s}} P^0_t f(s, \cdot)\, dt,
 \sL \varphi \right\> \,ds.
\end{split}
\end{equation}
By the Fubini theorem,
\begin{align*}
\E \int_0^{E_{T-s}} P^0_t f(s, \cdot)\,dt
&= \int_{t=0}^\infty \P(E_{T-s} >t) P^0_t f(s, \cdot) \,dt= \int_{t=0}^\infty \P(S_t\leq T-s) P^0_t f(s, \cdot) \,dt \\
&= \int_{t=0}^\infty \left( \int_{r=0}^{T-s} \bar p(t, r)\, dr\right)
P^0_t f(s, \cdot) \,dt= \int_{r=0}^{T-s} \int_{t=0}^\infty  P^0_t f(s, \cdot) \bar p(t, r) \,dt \,dr \\
&= \int_{r=0}^{T-s} \int_{t=0}^\infty  P^0_t f(s, \cdot) \bar p(t,
T-s-r)\, dt \,dr.
\end{align*}
Hence
\begin{align*}
\int_0^T \left( \E \int_0^{E_{T-s}} P^0_t f(s, \cdot)\,dt \right)
\,ds
&=\int_{s=0}^T \int_{r=0}^{T-s} \int_{t=0}^\infty  P^0_t f(s, \cdot) \bar p(t, T-s-r) \,dt\, dr\, ds \\
&=\int_{r=0}^T \int_{s=0}^{T-r} \int_{t=0}^\infty P^0_t f(s, \cdot) \bar p(t, T-r-s)\, dt\, ds \, dr \\
&= \int_0^T u(T-r, \cdot) \,dr.
\end{align*}
This combined with \eqref{e:2.8} shows that \eqref{e:2.7} holds for every $\varphi
\in \sD (\sL)$.
Therefore, $u(t, x)$ is a weak solution to
$\partial^w_t u = \sL u + f(t, x)$ on $(0,  T_0]\times E$.

 \medskip

 \noindent
{\rm (2)}
If \eqref{e:2.7n} holds, the conclusion (i') follows from Proposition \ref{P:2.3}(ii).
If, in addition, $\| \sL f(t, \cdot)\|_{L^2(E; \mu)} $
$\leq K<\infty$ for a.e. $t\in (0,  T_0]$,
the continuity of $t\mapsto \sL u(t,x)$ in $L^2(E; \mu)$ follows from
Proposition \ref{P:2.3}(iii).
\qed

\medskip

\noindent {\bf Proof of Theorem \ref{T:mainw}.}\,\,
That $u(t, \cdot)$ defined by \eqref{e:defu} is a weak solution of
\begin{equation}\label{e:2.28}
\partial^w_t u=
\sL u + f(t, \cdot) \quad \hbox{on $(0,  T_0]\times E$ with } u(0, \cdot)= g
\end{equation}
follows directly
from Proposition \ref{T:2.5}, \cite[Theorem 2.4]{CKKW} and the linearity.
So it remains to establish $u$ is the unique weak solution of \eqref{e:2.28}.

Suppose $v$ is another weak solution of \eqref{e:2.28}.
Then $h:=u-v$ is
a weak solution of
\begin{align}
\label{e:2.C21n}
\partial^w_t h = \sL h  \quad \hbox{on } [0,  T_0]\times E  \
\hbox{ with } \  h (0, \cdot )=0.
\end{align}
As in the proof of the uniqueness part of Theorem \ref{T:mains},
define $\bar f (t, x)$ and $\bar h (t, x)$ as
\eqref{e:dfbarf} and \eqref{e:dfbarh} respectively.
Since \eqref{e:2:C1}
holds in this case,
by  Proposition \ref{T:2.5}(1),
$\bar h(t, x)$ is a weak solution to
\begin{align}
\label{e:barhweak}
\partial^w_t \bar h = \sL \bar h (t, x)+ \bar f(t, x) \quad \text{on} \quad (0, \infty) \times E \quad \text{with} \quad \bar h(0, x)=0.
\end{align}
Moreover, \eqref{e:lambdabd} also holds.

Now define $h(t, x)= \bar h (t-T_0, x)$ for $t>T_0$.
Then,
by \eqref{e:TT0}, \eqref{e:2.C21n} and \eqref{e:barhweak},  we have that for all $T>T_0$ and  $\varphi \in \sD (\sL)$
\begin{align*}
 \int_0^T  \<\sL\varphi,    h (t, \cdot) \>\,dt
&=
\left\< \int_0^{T_0} w(T_0-s) h(t, \cdot)\,dt, \varphi \right\>
+
 \int_{T_0}^T  \<\sL\varphi,   \bar h (t-T_0, \cdot) \>\,dt\\
&=
\left\< \int_0^{T_0} w(T-s) h(t, \cdot)\,dt, \varphi \right\>+\left\<\int_{0}^{T_0}\left(w(T_0-t)-w(T-t)\right)h(t, \cdot)\,dt, \varphi \right\>\\
&\quad +
 \int_{0}^{T-T_0}  \<\sL\varphi,  \bar h (t, \cdot)\,dt \>\\
&=\left\< \int_0^{T_0} w(T-s) h(t, \cdot)\,dt, \varphi \right\> +\int_0^{T-T_0}\< \bar f(t, \cdot), \varphi \> \,dt + \int_0^{T-T_0} \<\bar h(s, \cdot), \sL
\varphi \> \,dt\\
 &= \left\< \int_0^{T_0} w(T-s) h(t, \cdot)\,dt, \varphi \right\>+\left\< \int_{T_0}^T w(T-s)  \bar h(t-T_0, \cdot)\,dt, \varphi \right\> \\
 &=\left\< \int_0^{T} w(T-s) h(t, \cdot)\,dt, \varphi \right\>.
\end{align*}
Hence, for every $t>0$ and $\varphi\in \sD (\sL)$,
\begin{equation}\label{e:2.3nn}
 \int_E \varphi(x) \left( \int_0^t w(t-r)  h(r, x) \,dr \right) \mu (dx)  = \int_E  h(t,x)\sL \varphi( x)   \,\mu (dx) .
\end{equation}

We now follow
the uniqueness proof in \cite[Theorem 2.3]{C} (see also \cite[Theorem 2.4]{CKKW})
to show that $h(t, x)=0$ in $L^2(E; \mu)$ for every $t>0$.
Let $H(\lambda, x):=\int_0^\infty e^{-\lambda t} h(t, x)\, dt$, $\lambda >0$,
 be the Laplace transform of $t\mapsto h(t, x)$. Note that, for every $\lambda >0$, $H(\lambda, \cdot)\in \bB$
 with $\| H(\lambda, \cdot)\| <\infty$; see \eqref{e:lambdabd}.
 Taking the Laplace transform in $t$ on both sides of \eqref{e:2.3nn} yields that for every $\lambda >0$,
\begin{align*}
\frac{\phi (\lambda)}{\lambda}\int_E \varphi(x) H(\lambda, x)\,\mu(dx)=&
\int_E \varphi(x) H(\lambda, x) \left(   \int_0^\infty e^{-\lambda s} w(s)\, ds \right) \,\mu(dx)\\
=& \frac{1}{\lambda} \int_E H(\lambda, x) \sL \varphi(x) \,\mu (dx).
\end{align*}
That is,  for every $\lambda >0$,
 $$
 \int_E H(\lambda, x)  \left( \phi (\lambda ) -\sL \right) \varphi(x) \,\mu (dx) =0.
 $$
 Denote by $\{G^0_\alpha: \alpha>0\}$ be the resolvent of the regular Dirichlet for $(\sE, \sF)$.
 For each fixed $\lambda >0$ and
$ \eta \in L^2(E; \mu)$, take $\varphi:=G^0_{\phi (\lambda)} \eta$,
which is
 in $\sD (\sL)$. Since $ ( \phi (\lambda ) -\sL ) \varphi  =\eta$, we deduce
 that $\int_E     H(\lambda , x )
 \eta (x) \,\mu (dx)=0$ for every $\eta\in L^2(E; \mu)$.
 Therefore $H(\lambda, \cdot )=0$ $\mu$-a.e. on $E$ for every $\lambda >0$. By the uniqueness of the Laplace transform
 and the fact that $t\mapsto h(t, x)$ is continuous in $L^2(E; \mu)$,
 it follows that $(u-v)(t, \cdot) = h(t, \cdot )=0$ in $L^2(E; \mu)$ for every $t>0$.
Therefore, we conclude that  $u(t, x)$ is the unique
 weak solution to
$\partial^w_t u = \sL u + f(t, x)$ on $(0,  T_0]\times E$ with $u(0, x)=g(x)$.
\qed

\section{Fundamental
solutions to time fractional Poisson equations}\label{section3}

In this section, we assume that \emph{the Banach space $(\bB, \| \cdot \|)$ is
a function space over a locally compact separable metric space $E$,
    where either $\bB= L^p(E; \mu)$ for some $\sigma$-finite measure $\mu$
    with full support on $E$ and $p\geq 1$, or $\bB = C_\infty (E)$,
    the space of continuous functions on $E$ that vanish at infinity.}
Let $\{P^0_t; t\geq 0\}$ be a
uniformly bounded  and strongly
    continuous semigroup in $(\bB, \| \cdot \|)$.
Recall
that $S=\{S_t, \P; t\geq 0\}$ is a driftless subordinator with infinite
L\'evy measure $\nu$ such that \emph{the subordinator $S_t$
  has a bounded density function $\bar p(t, \cdot )$ for each $t>0$}.
    As a particular case of
    Theorems \ref{T:mains} and \ref{T:mainw},
    we know that
    \begin{equation}\label{e:3.1}
    u(t, x) := \int_{s=0}^t \int_{r=0}^\infty P^0_r f (s, \cdot ) (x) \bar p (r, t-s)\, dr\, ds
    \end{equation}
    is the unique solution to the general time fractional  Poisson
        equation \eqref{e:3.2} with $u(0, x)=0$
     under suitable conditions on $f(t, x)$.
     We further assume the following assumption throughout this section.
\begin{assumption}\label{S:3.1}
 The uniformly bounded  and strongly
    continuous semigroup $\{P^0_t; t\geq 0\}$
 in $(\bB, \| \cdot \|)$
  has a density function $p^0(t, x, y)$ with respect to
 some $\sigma$-finite measure $\mu$ on $E$ with full support, such that
    \begin{description}
    \item{$(i)$} $P^0_t f(x) = \int_E p^0(t, x, y) f(y) \,\mu (dy)$ for $f\in \bB$;
    \item{$(ii)$} for each $x, y\in E$, $t\mapsto p^0(t, x, y)$ is Borel measurable.
    \end{description}
\end{assumption}

\medskip
    Under Assumption \ref{S:3.1}, $u(t, x)$ of \eqref{e:3.1} can be written as
    $$
    u(t, x)= \int_0^t \int_E q(t-s, x, y) f(s, y) \,\mu (dy) \,ds,
    $$
 where for $t>0$,
   \begin{equation}\label{e:3.3}
    q(t, x, y):= \int_0^\infty p^0(r, x, y) \bar p(r, t) \,dr.
   \end{equation}
    In other words, $q(t, x, y)$ is the fundamental
    solution  for solving the time fractional Poisson equation
 \eqref{e:3.2}
with zero initial value.

We know from
\cite[Theorem 2.3]{C}
(see also \cite[Theorem 2.4]{CKKW})
that
\begin{equation}\label{e:3.4}
p(t, x, y):= \E \left[ p^0 (E_t, x, y) \right]
\end{equation}
is the fundamental solution of the homogenous time fractional equation \eqref{e:3.2-0}.

Define $p(0, x, y)=0$ for $x\not=y$. This definition is reasonable in view of \eqref{e:3.4}
as $E_0=0$ and
$p^0(0, x, y)\,\mu(dy)=\delta_{\{x\}}(dy).$

Recall that  $G^{(S)} (t)$ is the potential density of the subordinator defined by \eqref{e:2.2a}.
\begin{prop}\label{P:3new}
 Suppose that $S_r$ has a
 bounded density  $\bar p(r, \cdot)$ for each $r>0$
and that Assumption $\ref{S:3.1}$ holds.
Then, for every $t>0$ and  $x, y\in E$,
\begin{equation}\label{e:3.5}
\int_0^t   q(s, x, y) \, ds = \int_0^t G^{(S)}(t-s) p(s, x, y)\, ds,
\end{equation}
\end{prop}
\pf
First note that by  Lemma \ref{T:2.1}, \eqref{e:3.3} and  \eqref{e:3.4},
 there is a zero Lebesgue set $\sN_0\subset (0, \infty)$ so that for every $t\in (0, \infty) \setminus \sN_0$ and
for every $x, y\in E$,
\begin{equation}\label{e:3.4b}
p(t, x, y) =\int_0^t w(t-s) q(s, x, y) \,ds .
\end{equation}
Thus, by  \eqref{e:3.4b}, the Fubini theorem and Lemma \ref{L:1.20}, for every $t>0$,
\begin{align*}
\int_0^t G^{(S)}(t-s) p(s, x, y) \,ds
&= \int_{s=0}^t G^{(S)} (t-s)  \left( \int_{r=0}^s w (s-r) q(r, x, y)\, dr \right)\, ds   \\
&= \int_{r=0}^t   q(r, x, y)  \left( \int_{s=r}^t w (s-r)  G^{(S)} (t-s) \,ds\right)\, dr    \\
&= \int_{r=0}^t   q(r, x, y)  \left( \int_{0}^{t-r} w (s)  G^{(S)} (t-r-s)\, ds\right)\, dr \\
&= \int_{ 0}^t   q(r, x, y)  \,dr.
\end{align*}
\qed

\medskip

Recall that the subordinator $S$ is said to be special, if its
Laplace exponent $\phi$ is a special Bernstein function; that is,
$\lambda\mapsto \phi^*(\lambda):= \lambda/\phi (\lambda)$ is still a
Bernstein function. In this case, let $S^*=\{S^*_t; t\geq 0\}$ be
the subordinator with the Laplace exponent $\phi^*(\lambda)$. We
call $S^*$ the conjugate subordinator to $S$. Let $\nu^*$ be the
L\'evy measure of $S^*$ and $w^*(x):= \nu^* (x, \infty)$.

The aim of this section is to show that, when $S$ is a special
subordinator and $S_r$ has a bounded density function $\bar p(r,
\cdot )$ for each $r>0$,  the fundamental solution
$p(t, x, y)$   of the homogenous time fractional equation \eqref{e:3.2-0}
has time-fractional
derivative $\partial^{w^*}_t p (\cdot , x, y) (t)$, and
$$
q(t, x, y) = \partial^{w^*}_t p (\cdot , x, y) (t)
\quad \hbox{for a.e. } t>0 \hbox{ and $\mu$-a.e. } x\not= y \hbox{ in } E.
$$
The following is a precise statement.

\begin{thm}\label{T:3.2} Suppose that $S$ is a special subordinator, that $S_r$ has a
 bounded density  $\bar p(r, \cdot)$ for each $r>0$
and that Assumption $\ref{S:3.1}$ holds.
 Denote by $\nu^*$ the L\'evy measure of the conjugate subordinator $S^*$ to $S$ and set
$w^*(x):=\nu^*(x, \infty)$.
Let $q(t,x,y)$ and $p(t,x,y)$ be defined  as in
\eqref{e:3.3} and \eqref{e:3.4}, respectively.  Then
\begin{equation}\label{e:3.4a}
\int_0^t q(s, x, y) \,ds = \int_0^t w^*(t-s) p(s, x, y)\, ds \quad \hbox{for } t> 0,
\end{equation}
which is finite for $\mu$-a.e. $x\in E$ and $\mu$-a.e. $y\in E\setminus \{x\}$.
Hence for any such $x\not=y$ so that \eqref{e:3.4a} holds,
$$
 \partial^{w^*}_t p(\cdot, x, y) (t ):=\frac{d}{dt} \int_0^t w^*(t-s) p(s, x, y) \,ds
\quad \hbox{exists for a.e. } t>0
$$
and
$$
q(t, x, y)=
\partial^{w^*}_t p(\cdot, x, y) (t)  \quad \hbox{for a.e. } t>0.
$$
\end{thm}

\pf
Since $S$ is a special subordinator and $\nu$ is infinite, $w^*(t)=G^{(S)}(t)$
(see \cite[Corollary 11.8]{SSV}). Thus \eqref{e:3.4a} follows immediately from this and Proposition \ref{P:3new}.

We now show the finiteness of \eqref{e:3.4a}.
Let $M=\sup_{t>0}\|P^0_t\|$. For every $t>0$ and $f\in \bB$, by
\eqref{e:3.1},
 \eqref{e:3.3} and Fubini's theorem,
\begin{align*}
 \left\| \int_E   \left(\int_0^tq(s, \cdot , y)\,ds\right) f(y)\,
   \mu(dy) \right\|
&\leq   \int_0^t  \int_0^\infty \| P^0_r f \| \, \bar p(r, s) \,dr\, ds  \\
&\leq M \| f\| \int_0^t  \int_0^\infty \bar p(r, s) \,dr \,ds
  =M \| f\|  G^{(S)}{\bf1}_{(0, t]}(0)<\infty.
\end{align*}
Hence by \eqref{e:3.4a}, for $\mu$-a.e. $x\in E$ and $\mu$-a.e. $y\in E\setminus \{x\}$,
$$
\int_0^t q(s, x, y) \,ds =\int_0^t w^*(t-s) p(s, x, y)\, ds <\infty
\quad \hbox{for every } t>0.
$$
This proves the theorem.
 \qed

\section{Two-sided estimates}\label{section4}

\subsection{Density estimates for subordinators}
In this subsection, we assume that the  driftless subordinator
$S=\{S_t, \P;t\ge0\}$ satisfies the following
assumptions.

\begin{assumption}\label{as:sub}
 \rm
{\rm (i)} The Laplace exponent $\phi$ of $S$ satisfies that
\begin{equation}\label{e:phi}
c_1 \kappa^{\beta_1}\le \frac{\phi(\kappa \lambda)}{\phi(\lambda)}\le c_2 \kappa^{\beta_2}
\quad\text{ for all }\lambda>0 \text{ and }\kappa\ge 1,
\end{equation}
where $0<\beta_1<\beta_2<1$. Without loss of generality, we assume
$\phi(1)=1$.\\
{\rm (ii)}  The L\'evy measure $\nu(dz)$ of $S$ has a density function $\nu(z)$ with respect to the Lebesgue measure such that
 \begin{align}
 \label{e:tvdec}
 t \mapsto t\nu(t) \text{ is   non-increasing
  on } (0, \infty).
 \end{align}
\end{assumption}

We make some remarks on the assumption above.

\begin{remark}\label{remark4.2} \rm \begin{description}
\item{(i)} Under \eqref{e:phi}, the L\'evy measure $\nu$ of $S$ is infinite as
  $\nu(0,\infty)=\lim\limits_{\lambda \to \infty}\phi(\lambda)=\infty$,
 excluding compound Poisson processes.

\item{(ii)} Condition \eqref{e:tvdec}
is equivalent to saying that $\P(S_1 \in dt)$ is self-decomposable (see \cite[Definition 5.14 and Proposition 5.17]{SSV}).

\item{(iii)}
Under \eqref{e:phi} and \eqref{e:tvdec},
it follows from the proof of  \cite[Proposition 2.5]{KSV}
that
$\nu(t)
\simeq
 t^{-1}\phi(t^{-1})$ for all $t>0.$

\item{(iv)} According to \cite[Theorems 27.13 and 28.4(ii)]{Sat},
\eqref{e:phi} and \eqref{e:tvdec} together imply that
 $S_r$ has
a bounded
 density $\bar p(r, t)$ so that $\P(S_r \in dt)=\bar p(r, t)\, dt$; moreover, $t \mapsto \bar p(r, t)$ is smooth for any $r>0$.
\end{description}
\end{remark}

We now mention two large classes of subordinators that satisfy both \eqref{e:phi} and \eqref{e:tvdec}.

 \begin{example} \rm (1)
Let $0<\beta_1<\beta_2<1$. Suppose  a positive function
$\kappa(\beta, t)$ defined on $[\beta_1, \beta_2] \times (0,
\infty)$ satisfies that
 $c_1^{-1} \le \kappa(\beta, t)\le c_1$ for all $(\beta,t)\in [\beta_1,\beta_2]\times(0,\infty)$ and some constant $c_1\ge1$, and that
 $t \mapsto \kappa(\beta, t)$  is
 non-increasing
 on $(0, \infty)$ for any fixed $\beta\in [\beta_1,\beta_2]$. Define
 $$\phi(\lambda):=\int_{0}^\infty (1-e^{-\lambda t} )\nu(t)\,dt,$$
 where $$
 \nu(t):=t^{-1} \int_{\beta_1}^{\beta_2}\kappa(\beta, t)
 t^{-\beta}\,\mu_I(d\beta)$$ and $\mu_I$ is a finite measure on
 $I:=[\beta_1,\beta_2]$.
 Then clearly \eqref{e:tvdec} holds.
 Furthermore, since
 $$\phi(\lambda)
\simeq \int_{\beta_1}^{\beta_2}
 \int_{0}^\infty (1-e^{-\lambda t} )t^{-1-\beta} \,dt\,\mu_I(d\beta)
 =\int_{\beta_1}^{\beta_2} \left(\int_0^\infty(1-e^{-s})s ^{-1-\beta}\,ds \right)
 \lambda^\beta\,\mu_I(d\beta),
 $$
 it is easy to see that  \eqref{e:phi} holds too.
In particular, stable subordinators, which correspond to $\mu_I$ being a Dirac measure $\delta_{\{\beta\}}$
for some $\beta\in(0,1)$ and
$\kappa (\beta, t) \equiv c_0$ for all $t>0$ for some positive constant $c_0$,
and a larger class of mixtures of stable subordinators satisfy both
 \eqref{e:phi} and \eqref{e:tvdec}. For more related examples, see \cite[Examples 15.13--15.15]{Sat} and the references
 therein.

(2) A function $f:(0,\infty)\to [0,\infty)$ is said to be completely monotone,  if $f$ is of class
$C^\infty$ and $(-1)^nf^{(n)}\ge 0$ on $(0,\infty)$ for every integer $n\geq 0$.
A Bernstein function is said to be a complete Bernstein function, if its
L\'evy measure has a completely monotone density  with respect to
the Lebesgue measure.
A sufficient condition on $\phi$ which implies
\eqref{e:tvdec} is that $\phi$ is a Thorin-Bernstein function; that
is, both $\phi(\lambda)$
and $\lambda\phi'(\lambda)$ are complete
Bernstein functions. In this case,  both $t\mapsto \nu (t)$ and $t \mapsto t \nu (t)$ are completely
monotone. (See \cite[Definition 8.1, Theorem
8.2(iv) and Theorem 6.2]{SSV}.)
\end{example}

The main result of this subsection is as follows.
\begin{thm}\label{t:sub} Let $S$ be a subordinator such that
\eqref{e:phi}
and \eqref{e:tvdec} are satisfied. Denote by $\bar p(r,t)$
the associated transition density.
Then, for each $L>0$, there exist
constants $c_i:=c_{i,L}\ge1$ $(i=1,2,3)$ such that
\begin{equation}\label{e:hksubt1}
\frac{1}{c_1t} r\phi(t^{-1}) \le \bar p(r,t)\le \frac{c_1}{t}
r\phi(t^{-1}) \quad \text{ for every positive } r, t \text{ with
}r\phi(t^{-1})\le L
\end{equation}
and
\begin{equation}\label{e:hksubt2}
\frac{1}{c_2 t}e^{-c_3t(\phi')^{-1}(t/r)} \le \bar p(r,t) \le
\frac{c_2}{t}e^{-c_3^{-1}t(\phi')^{-1}(t/r)} \quad \text{ for every
positive } r, t \text{ with }  r\phi(t^{-1}) \ge  L,
\end{equation}where $(\phi')^{-1}(\lambda):=\inf\{s>0: \phi'(s)\le \lambda\}$ for all $\lambda\ge0.$
\end{thm}

 \ \

To prove Theorem \ref{t:sub}, we will make use of the following scaled processes. For $r>0$, define
$$
\phi^r(\lambda):=\frac{\phi(\lambda r^{-1})}{\phi(r^{-1})}\, ,\quad \lambda >0\, .
$$
Note that $\phi^r(1)=1$ for all $r>0$. Since
$$\phi^r(\lambda)=\frac{\phi(\lambda r^{-1})}{\phi(r^{-1})}= \int_{0}^\infty (1-e^{-\lambda t/r} )\phi(r^{-1})^{-1}\nu(t)\,dt=\int_{0}^\infty (1-e^{-\lambda s} )\phi(r^{-1})^{-1}r\nu(rs)\,ds,
$$
$\phi^{r}$ is a driftless Bernstein function with  L\'evy density $\nu^r$ given by
$ \nu^r(t)=r\phi(r^{-1})^{-1}\nu(r t)$.
It is obvious that \eqref{e:tvdec} is satisfied, and \eqref{e:phi} holds for all $\phi^{r}$ with the same constants, i.e.,
\begin{equation}\label{e:phic} c_1 \kappa^{\beta_1}\le \frac{\phi^r(\kappa \lambda)}{\phi^r (\lambda)}\le c_2 \kappa^{\beta_2}\quad \text{for any  }\lambda, r >0 \text{ and } \kappa\ge 1.
\end{equation}
In particular, according to \cite[Proposition 2.5]{KSV} (and its proof) again, \eqref{e:phic} implies that
\begin{equation}\label{e:nuc2} \nu^r(t)
\simeq
\phi^r(t^{-1})t^{-1}, \quad t, r>0.
\end{equation}

Let $S^r=\{S^r_t, \P;t\ge 0\}$ be a subordinator with the Laplace exponent $\phi^r$.
Since
$$
\E[e^{-\lambda S^r_t}] =e^{-\frac{t}{\phi(r^{-1})} \phi(\lambda r^{-1}) }  =\E[e^{-\lambda r^{-1} S_{t/ \phi(r^{-1})}}],
$$ $S^r$ is identical in law to the process
$\{r^{-1}S_{t/ \phi(r^{-1})}, \P; t\ge0 \}$.
In particular,
$ \frac{1}{ \phi^{-1}(1/t) }S_1^{1/ \phi^{-1}(1/t) }$ enjoys the same law as that of $ S_{t}$, so
the density $q^r$ of  $S_1^r$ has the identity
\begin{align}
\label{e:qvbp}
\phi^{-1}(1/t) q^{1/ \phi^{-1}(1/t) }(u\phi^{-1}(1/t))=\bar p(t, u), \quad t,r>0.
\end{align}

We have the following strictly unimodal property and
an upper bound for $\bar p(r,t).$

\begin{lemma}
\label{t:spd} Under \eqref{e:phi},
and  \eqref{e:tvdec}, the following hold.
\begin{itemize}
\item[$(i)$]
There exists a constant $a_r>0$
such that for all $r>0$,  $t\mapsto \bar p(r, t)$ is strictly increasing on $[0,a_r)$ and $t\mapsto \bar p(r, t)$ is strictly decreasing on $(a_r, \infty)$.
Moreover, there exists a constant $c_1\ge1$ such that
$c_1^{-1}/\phi^{-1}(1/r) \le a_r \le  c_1/\phi^{-1}(1/r)$ for all $r>0$.
\item[$(ii)$]  There exists a constant $c_2>0$ such that for all $r,t>0$,
$$\bar p(r,t)   \le c_2  \phi^{-1}(1/r).$$
\end{itemize}
\end{lemma}
\pf (i) We first note that, since $\lim_{t \downarrow 0} t \nu(t)=\infty$ and $\int_{0}^1 t \nu(t)\,dt<\infty$ by \eqref{e:phi},  \eqref{e:tvdec} and Remark \ref{remark4.2}(iii),
 $S_r$ is of type $I_6$ in \cite{SY} for each $r>0$
(see \cite[p.\ 275]{SY}). Thus,
according to \cite[Theorem 1.4]{SY} and \cite[Corollary 1]{W},  $\P(S_r \in dt)$ is strictly unimodal with mode $a_r \ge 0$ on $[0, \infty)$; that is,
$\P(S_r \in dt)=c_r \delta_{a_r} (dt) + \wh p(r, t) \,dt$,
where $c_r \ge 0$, $t\mapsto \wh p(r, t)$ is strictly increasing on $[0,a_r]$ and $t\mapsto \wh p(r, t)$ is strictly decreasing on $[a_r, \infty)$.
Since $\P(S_r \in dt)=\bar p(r, t) \, dt$, we have $c_r = 0$ and so $\wh p(r, t)=\bar p(r, t)$.

We next show that $a_r \asymp 1/\phi^{-1}(1/r)$ on $(0, \infty)$ using the above scaled processes $S^{r}$.
We first claim that the  infinitely
 divisible random variables $S_1^{r}$'s have their modes $b_r$'s such that $b_r
\simeq 1$ for all $r>0$. First of all,
since
$ t \mapsto t\nu^r(t)=r\phi(r^{-1})^{-1}t\nu(r t)$ is  non-increasing on $(0, \infty)$,  $S_1^{r}$'s are self-decomposable by  Remark \ref{remark4.2}(ii);
moreover, $S_1^{r}$'s are of type $I_6$ in \cite{SY} because of \eqref{e:phic}, \eqref{e:nuc2} and Remark \ref{remark4.2}(iii).
Thus, according to \cite[Theorem 6.1 (ii) and (vi)]{SY}, we have that
$b_r^{-1}\int_0^{b_r} \nu^r(t)t \,dt >1$ and $b_r >\sup\{t: \nu^r(t)t \ge 1\}$.
Using $\phi^r(1)=1$, \eqref{e:phic} and \eqref{e:nuc2} to these two inequalities, we get
$$b_r >\sup\{u: \nu^r(u)u \ge 1\} \ge c_1 \sup\{u:  \phi^r(u^{-1}) \ge  1\} \ge c_1$$
and
\begin{align*}
1&<b_r^{-1}\int_0^{b_r} \nu^r(t)t \,dt  \le c_2 b_r^{-1}\int_0^{b_r} \phi^r(t^{-1}) \,dt
 = c_2   \phi^r(b_r^{-1}) b_r^{-1}\int_0^{b_r} \frac{\phi^r(t^{-1})} {\phi^r(b_r^{-1})} \,dt \\
 & \le c_3  \phi^r(b_r^{-1}) b_r^{-1}\int_0^{b_r} \left(\frac {b_r}{t} \right)^{\beta_2}\,dt
 = c_4 \phi^r(b_r^{-1}).
\end{align*}
Therefore, by the inequalities
above and \eqref{e:phic},
$$
c_1 \le b_r \le \frac{1}{\phi_r^{-1}(1/c_4)}=\frac{\phi_r^{-1}(1)}{\phi_r^{-1}(1/c_4)} \le c_5.
$$
We have proved the claim.
Now, using \eqref{e:qvbp}, we conclude that $a_r
\simeq 1/\phi^{-1}(1/r)$ on $(0, \infty)$.

(ii) Recall that  $q^r$ is the density of $S_1^r$.
Let $$x^r(\lambda):=\int_0^\infty (1-\cos(\lambda t))\nu^r(t)\,dt$$ and
$$y^r(\lambda):=\int_0^\infty \sin(\lambda t)\nu^r(t)\,dt,$$ which are finite by  \eqref{e:phic} and \eqref{e:nuc2}.
\eqref{e:phic} and \eqref{e:nuc2} also imply that the characteristic exponent of $S^r_t$ is
$$\psi^r(\lambda):=x^r(\lambda)-iy^r(\lambda)=\int_0^\infty (1-e^{i\lambda t})\nu^r(t)\,dt, \quad \lambda>0.$$
Using the inequality $(1-\cos s) \ge (1-\cos 1)s^2$ for any $s \in (0, 1)$, \eqref{e:phic} and
\eqref{e:nuc2},  we have
\begin{align*}
x^r(\lambda) &\ge (1-\cos 1) \int_0^{1/\lambda} (\lambda t)^2\nu^r(t)\,dt
\ge c_1 \lambda^2 \phi^r(\lambda) \int_0^{1/\lambda} t \frac{\phi^r(1/t)}{\phi^r(\lambda)}\,dt \\
&\ge c_2 \lambda^2 \phi^r(\lambda) \int_0^{1/\lambda} t (t\lambda)^{-\beta_1}\,dt  \ge c_3 \phi^r(\lambda).
\end{align*}
In particular,  $e^{-\psi^r} \in L^1(\bR;dx)$,  so by applying the Fourier inversion formula and then using \eqref{e:phic} we can bound $q^r$ as
\begin{align*}
q^r(z)=&\frac{1}{2\pi}\int_{\bR} e^{-\lambda z-\psi^r(z)}\,dz
=\frac{1}{2\pi}\int_{\bR} e^{i(y^r(\lambda)-\lambda z) - x^r(\lambda)   } \,d\lambda\\
 \le &\frac{1}{\pi} \int_0^\infty
e^{-x^r(\lambda) }\,d\lambda
\le  \frac{1}{2\pi} \left( 1+ \int_1^\infty
e^{-c_4 \lambda^{\beta_1} } \,d\lambda \right).
\end{align*}
Thus, $q^r$'s are uniformly bounded. Now,
the second assertion follows from this and \eqref{e:qvbp}.\
\qed

To establish two-sided estimates of $\bar p(r, t)$, we further need the following two statements from \cite[Lemma 3.1 and Proposition 3.3]{CKKW}, respectively.

\begin{lemma}[{\bf\cite[Lemma 3.1]{CKKW}}]
\label{l:AKM} Let $\phi$ be a Bernstein function such that
\eqref{e:phi} is satisfied. Then there exists a constant $C_*\ge 1$
such that the following holds
\begin{equation}\label{e:04}
 \lambda \, \phi'(\lambda) \le \phi(\lambda) \le C_* \lambda \, \phi'(\lambda),
\quad \lambda>0.
\end{equation}
In particular, there exist constants $c_i>0$ $(i=3,4,5,6)$ such that
\begin{equation}\label{e:031}c_3 \kappa^{1-\beta_2}\le \frac{\phi'(\lambda)}{\phi'(\kappa
\lambda)}\le c_4 \kappa^{1-\beta_1},\quad \lambda>0, \kappa\ge
1,\end{equation} and
\begin{equation}\label{e:03}c_5 \kappa^{1/(1-\beta_1)}\le \frac{(\phi')^{-1}(\lambda)}{(\phi')^{-1}(\kappa \lambda)}\le c_6 \kappa^{1/(1-\beta_2)},\quad \lambda>0, \kappa\ge 1.\end{equation}
\end{lemma}

\begin{prop}[{\bf \cite[Proposition 3.3]{CKKW}}]\label{p:sub}  Let $S=\{S_t,\P;t\ge0\}$ be a subordinator whose Laplace exponent
 $\phi$ satisfies assumption \eqref{e:phi}. Then, we have
\begin{itemize}
\item[$(i)$] There are constants $c_1,c_2>0$ such that for all $r,t\ge0$,
\begin{equation*}
\bP  \left( S_r\ge t(1+er\phi(t^{-1}) \right) \le c_1 r\phi(t^{-1})
\end{equation*}
and
\begin{equation*}
\bP (S_r\geq  t)
\ge 1- e^{-c_2r \phi(t^{-1})}.\end{equation*}
 In particular, for each $L>0$, there exist constants $c_{1,L},c_{2,L} >0$ such that
for all $r\phi(t^{-1})\le L$,
$$
c_{1,L} r\phi(t^{-1})\le \bP (S_r\ge t)\le  c_{2,L} r\phi(t^{-1}).
$$

\item[$(ii)$] There is a constant $c_1>0$ such that for all $r, t>0$,
$$ \bP (S_r\le t)\le  \exp(-c_1 r \phi\circ[(\phi')^{-1}] (t/r)) \le  \exp(-c_1t (\phi')^{-1}(t/r)).$$
Moreover, there is a constant $c_0>0$ such that for each $L>0$, there exists a constant $c_{c_0,L}>0$ so that  for  $r\phi(t^{-1}) > L$
 $$ \bP (S_r\le t)\ge c_{c_0,L} \exp\left(-c_0 r \phi\circ[(\phi')^{-1}] (t/r)\right)\ge  c_{c_0,L} \exp(- c_0 C_*t (\phi')^{-1}(t/r)),$$ where $C_*>0$ is an absolute constant only depending on $\phi$.\end{itemize}
\end{prop}

Now, we are in a position to present the

 \ \

\noindent {\bf Proof of Theorem \ref{t:sub}.}\,\, (1)
By Lemma \ref{t:spd}(i), there exists a constant $L\ge1$ such that  for every fixed $r>0$,
 $t\mapsto\bar p(r,t)$ is strictly decreasing for $t \ge 1/ \phi^{-1}(1/(Lr))$ and
  $t\mapsto\bar p(r,t)$ is strictly increasing for $t \le 1/ \phi^{-1}(L/r)$.

 By Proposition \ref{p:sub}(i), for any $r,t>0$ with $r\phi(t^{-1})\le 1/ L$ and for any $c>1$,
\begin{align*}
\bar p(r,t)  \ge & \frac{1}{(c-1)t}\int_t^{ct}\bar p(r,s)\,ds=\frac{1}{(c-1)t}\left(\bP(S_r>t)-\bP(S_r>ct)\right)\\
\ge&\frac{r}{(c-1)t}\left(c_{1,L} \phi(t^{-1})-c_{2,L} \phi((ct)^{-1})\right).\end{align*}
In the first inequality above we have used the fact that $s\mapsto\bar p(r,s)$ is
decreasing
 on $(t, ct)$.
Thus, by \eqref{e:phi}, taking $c>0$ large enough, we arrive at
$$ \bar p(r,t)\ge
\frac{c_{1,L}}{2(c-1)t}
r\phi(t^{-1}) \quad \text{for all } r\phi(t^{-1})\le 1/L.$$
Similarly,
since  $s\mapsto\bar p(r,s)$ is
decreasing
 on $(t/2,\infty)$ if $r\phi( 2t^{-1})\le 1/L$, we have that
for  $r\phi( 2t^{-1})\le 1/L$
\begin{align*}
\bar p(r,t)  \le \frac{2}{t}\int_{t/2}^{t}\bar p(r,s)\,ds=\frac{2}{t}\left(\bP(S_r>t)-\bP(S_r>t/2)\right)\le\frac{2}{t} c_{2,L} r\phi(t^{-1}).\end{align*}

On the other hand, according to Proposition \ref{p:sub}(ii), for any $r,t>0$ with $r\phi(t^{-1})> L$
\begin{align*}
\bar p(r, t) \ge \frac{1}{t}\int_{0}^{t}\bar p(r,s)\,ds=\frac{1}{t}\bP(S_r\le t) \ge\frac{1}{t}c_{c_3,L} \exp(- c_{4}t (\phi')^{-1}(t/r)).\end{align*}
In the first inequality above we have used the fact that $s\mapsto\bar p(r,s)$ is increasing on $(0, t)$.
Similarly,
since  $s\mapsto\bar p(r,s)$ is increasing  on $(0,2t)$ if $r\phi((2t)^{-1})> L$, we have that
for any $r,t>0$ with $r\phi((2t)^{-1})> L$,
\begin{align*}
\bar p(r, t) \le \frac{1}{t}\int_t^{2t}\bar p(r,s)\,ds=\frac{1}{t}\left(\bP(S_r\le 2t)-\bP(S_r\le t)\right) \le \frac{1}{t} \exp(-c_5t (\phi')^{-1}(t/r)).
\end{align*}
Hence, combining above estimates and using \eqref{e:phi}, we have that there exists a constant $C>1$ such that
\begin{equation}\label{e:hksub0}
\bar p(r,t)
\simeq \frac{1}{t} r\phi(t^{-1}) \quad \text{ for every positive } r, t \text{ with
}r\phi(t^{-1})\le C^{-1}
\end{equation}
and
\begin{equation}\label{e:hksub1}
\frac{c_{6,C}}{t}e^{-c_7 t(\phi')^{-1}(t/r)} \le \bar p(r,t) \le \frac{1}{t}e^{- c_8   t(\phi')^{-1}(t/r)} \quad \text{ for every positive } r, t \text{ with
} r\phi(t^{-1}) \ge  C.
\end{equation}

\noindent
(2) Next, we consider the case that $t,r>0$ satisfying $C^{-1} \le r\phi(t^{-1}) \le  C$.
Since
$$r \cdot \phi'\circ[(\phi')^{-1}](t/r) \cdot (\phi')^{-1}(t/r)= t \cdot (\phi')^{-1}(t/r),$$
by \eqref{e:04} we have
\begin{equation}\label{e:hksub3}
 t \cdot (\phi')^{-1}(t/r)
\simeq \frac{\phi\circ[(\phi')^{-1}](t\phi(t^{-1}))}{\phi(t^{-1})}=  \frac{\phi\circ[(\phi')^{-1}]\left(\frac{\phi(t^{-1})}{t^{-1}}\right)}{\phi\circ[(\phi')^{-1}] (\phi'(t^{-1}))}
\simeq 1. \end{equation}
Using Proposition  \ref{t:spd} and \eqref{e:hksub3},
we get
\begin{equation}\label{e:hksub4}\bar p(r,t) \le  c_{1}  \phi^{-1}(1/r) \le c_{1}  \phi^{-1}(C \phi(1/t)) \le c_{2}t^{-1} \le c_3
t^{-1} [C  r\phi(t^{-1}) \wedge e^{-c_{4}t(\phi')^{-1}(t/r)}].
\end{equation}
 For lower bound of $\bar p(r,t)$ on $C^{-1} \le r\phi(t^{-1}) \le  C$, we use
the unimodality of $\bar p(r,t)$ and \eqref{e:hksub0}--\eqref{e:hksub3}, and get
\begin{equation}\label{e:hksub5} \bar p(r,t) \ge \bar p(r, 1/ \phi^{-1}(C/r)) \wedge \bar p(r, 1/ \phi^{-1}(r/C))
\ge c_{5} t^{-1} \ge  \frac{c_{5}}{2 t}( C^{-1} r\phi(t^{-1})+ e^{-t(\phi')^{-1}(t/r)}).
\end{equation}

Therefore, combining \eqref{e:hksub0}, \eqref{e:hksub1}, \eqref{e:hksub4} and  \eqref{e:hksub5} together, we have proved the theorem.
 \qed

\subsection{Two-sided estimates
of the fundamental solution}
Suppose that $0<\alpha_1 \le \alpha_2 < \infty$.
We say that a non-decreasing function $\Psi: (0,\infty)\to (0,\infty)$ satisfies the \emph{weak scaling property with $(\alpha_1,\alpha_2)$}
if there exist constants $c_1$ and $c_2>0$
such that \begin{equation}\label{e:lleeqpa}
c_1(R/r)^{\alpha_1}
\le\Psi(R)/\Psi(r) \le c_2(R/r)^{\alpha_2}
\quad \text{for all } 0<r \le R < \infty.
\end{equation}
We say that a family of non-decreasing functions $\{\Psi_x\}_{x\in \Lambda}$  satisfies the \emph{weak scaling property uniformly with $(\alpha_1,\alpha_2)$} if each $\Psi_x$ satisfies the weak scaling property with constants
$c_1,c_2>0$ and $0<\alpha_1 \le \alpha_2 < \infty$ independent of the choice of $x\in \Lambda$.

In this subsection, let $( E,d, \mu ) $ be a locally compact separable metric
measure space such that $\mu$ is a Radon measure on $( E,d )$  with
full support.
 For $x\in E$ and $r\ge 0$, define
\[
V(x,r)=\mu(B(x,r)).
\]
We further assume that for each $x\in E$, $V(x,\cdot)$ satisfies the weak scaling property
uniformly with $(d_1,d_2)$ for some $d_2 \ge d_1>0$;
that is, for any $0<r\le R$ and $x\in E$,
\begin{equation}\label{vd1} c_1\left(\frac R r\right)^{d_1}\le \frac{V(x,R)}{V(x,r)}\le c_2 \left( \frac R r\right)^{d_2}.\end{equation}
Note that \eqref{vd1} is equivalent to the so-called volume doubling and reverse volume doubling
conditions.

Let $\{P^0_t; t\geq 0\}$ be  a
uniformly bounded and strongly continuous
semigroup on $L^p(E;\mu)$ with $p\ge 1$ or on $C_\infty(E)$
having a density  $p^0(t, x, y)$ with respect to $\mu$. Let $\{S_t, \P;t\ge0\}$ be the driftless subordinator
 satisfying \eqref{e:phi}
 and  \eqref{e:tvdec}. In particular, $S_r$ has a density $\bar p(r, t)$ so that
 $\bar p(r, t)$
 obeys two-sided estimates as these in Theorem \ref{t:sub}.

The purpose of this part is to establish
two-sided estimates for
 the fundamental solution $q(t, x, y)$ for the time fractional Poisson equation
 given by \eqref{e:3.3}, i.e.,
$$
q(t, x, y) = \int_0^\infty p^0(r,x,y) \bar p (r, t)\,dr.$$
We adopt the setting of \cite{CKKW}, where two-sided estimates on
the fundamental solution $p(r, x, y)$ of the homogenous time fractional equation
were
derived. We will consider estimates for $q(t,x,y)$ when $\{P^0_t; t\geq 0\}$ is associated with pure jump processes and diffusion processes, respectively.
Since the arguments are partly similar to these in \cite[Sections 4 and 5]{CKKW},  we just present main results on two-sided estimates for
the fundamental solution for the time fractional Poisson equation,
and then give some intuitive explanations on these estimates
and present the proof of Theorem \ref{T:1.4}.
Full proofs of theorems are postponed to the appendix of this paper.

\subsubsection{{\bf Pure jump case}}
We first consider the case that $\{P^0_t; t\geq 0\}$ is corresponding to
a strong Markov processes with pure jumps, where the associated transition
density $p^0(t, x, y)$ enjoys the following two-sided estimates:
\begin{equation}\label{e:hkmequi}p^0(t,x,y)\simeq
\frac{1}{V(x, \Phi^{-1}(t))} \wedge \frac{t}{V(x,
d(x,y))\Phi(d(x,y))}=:
{\overline p^0}(t,x,d(x,y)), \quad t>0, x,y\in E.
\end{equation} Here $\Phi :[0,+\infty )\rightarrow \lbrack
0,+\infty )$ is a strictly increasing function with $\Phi(0)=0$ that
satisfies the weak scaling property with $(\alpha_1,\alpha_2)$,
i.e., \eqref{e:lleeqpa} is satisfied.
The most typical examples are the heat kernels for $\alpha$-stable-like processes on $d$-sets, in which
$\Phi (r)=r^\alpha$ and $V(x,r)=r^d$ for all $x\in E$ and $r>0$. Further
examples of pure jump processes with heat kernel estimates can be referred to \cite{CK, CKW} and the references therein.

\begin{thm}\label{theorem:mainjump}
We assume that Assumption $\ref{as:sub}$ holds.
 Let $q(t,x,y)$ be given by
\eqref{e:3.3}.
Suppose that $p^0(t,x,y)$ enjoys
two-sided estimates \eqref{e:hkmequi}.
Then we have the following two statements:
\begin{itemize}
\item[$(i)$]
If   $\Phi(d(x,y))\phi(t^{-1}) \le 1$, then
\begin{equation}\label{e:llffoo}\begin{split} q(t,x,y) &\simeq  \frac{ \phi(t^{-1})}{t} \int_{\Phi(d(x,y))}^{2/\phi(t^{-1})} \frac{r}{V(x,\Phi^{-1}(r))}\,dr.
        \end{split}
        \end{equation}

\item[$(ii)$]
If   $\Phi(d(x,y))\phi(t^{-1}) \ge 1$, then $$q(t,x,y)\simeq
  \frac 1{t\phi(t^{-1})^2 V(x,d(x,y))\,  \Phi(d(x,y)) }.
 $$
\end{itemize}
\end{thm}

\noindent {\bf Proof of Theorem \ref{T:1.4}\,(i).}\,\,
By \cite[Theorem 3.2]{GHL}, our assumption on the estimates
\eqref{e:1.7}  of Markovian transition kernel $p^0(t, x, y)$
implies that $V(x,r)\simeq r^d$.
Thus, $\Phi(r)=r^\alpha, V(x,r)\simeq r^d$ and $\phi(r)=r^\beta$ in this case,
Therefore,
the condition $\Phi(d(x,y))\phi(t^{-1}) \le 1$ is equivalent to
$d(x,y)\le t^{\beta/\alpha}$. The conclusion of  Theorem \ref{T:1.4}\,(i) is
then an easy consequence of Theorem \ref{theorem:mainjump}. \qed

\subsubsection{{\bf Diffusion case}}

We next consider the case that
the strong Markov process $X$ corresponding to $\{P^0_t; t\geq 0\}$
is a conservative diffusion;
that is, it has continuous sample paths and infinite lifetime.
We assume that
the heat kernel $p^0(t,x,y)$
of the diffusion $X$ with respect to $\mu$ exists and  enjoys the
following two-sided estimates
\begin{equation}\label{eq:fibie3}
p^0(t,x,y)\asymp \frac{1}{V(x, \Phi^{-1}(t))}
\exp\left(-m(t,d(x,y))\right)=:{\overline p^0}(t, x,d(x,y)), \quad t>0, x,
y\in E.
\end{equation} Here, $\Phi :[0,+\infty )\rightarrow \lbrack
0,+\infty )$ is a strictly increasing function with $\Phi(0)=0$, and
satisfies the weak scaling property with $(\alpha_1,\alpha_2)$ such
that $\alpha_2 \ge \alpha_1>1$ in \eqref{e:lleeqpa};
the function $m(t,r)$ is strictly positive for all $t,r>0$,
  non-increasing  on $ t \in (0,\infty)$
 for fixed $r>0$, and  determined by
\begin{equation}\label{e:scdf}
\frac{t}{m(t,r)}\simeq \Phi\left(\frac{r}{m(t,r)}\right),\quad t,r>0.
\end{equation}
In all the literature we know,
for example
\cite[Section 7]{Grig}, \cite[Page 1217--1218]{GT} and \cite{HK}),
estimate \eqref{eq:fibie3}
was established under assumptions that include $( E,d ) $ being connected and satisfying the chain condition
(that is,
there exists a constant $C>0$ such that, for any $x,y\in
M$ and for any $n\in {\mathbb N}$, there exists a sequence
$\{x_{i}\}_{i=0}^{n}\subset E$ such that $x_{0}=x$, $x_{n}=y$ and
$
d(x_{i},x_{i+1})\leq C { d(x,y)}/{n}$ for all $i=0,1,\cdots,n-1$).
However, given the estimate \eqref{eq:fibie3} on $p_0(t, x, y)$,
we do not need to assume the connectedness nor the chain condition on $(E, d)$
in arguments below
for  the estimates of $q(t, x, y)$.
Since \eqref{e:lleeqpa} holds with $\alpha_1>1$ and \eqref{e:scdf} is satisfied, there are constants $c_1,c_2>0$ such that for all $r>0$,
\begin{equation}\label{e:scdf-1}c_1\left(\frac{T}{t}\right)^{-1/(\alpha_1-1)}\le \frac{m(T,r)}{m(t,r)}\le c_2\left(\frac{T}{t}\right)^{-1/(\alpha_2-1)},\quad 0<t\le T.
\end{equation}
On the other hand, due to \eqref{e:scdf},
 \begin{align}
 \label{e:mphisim1}
 m(\Phi(r),r)\simeq1, \quad  r>0.
 \end{align}
Using this and the fact that  $m(\cdot,r)$ is non-increasing,
we have that for every $c_3>0$ there exists $c_4\ge 1$ such that
$$
\frac{c_4^{-1}}{V(x, \Phi^{-1}(t))} \le p^0(t,x,y)\le  \frac{c_4}{V(x, \Phi^{-1}(t))}
\qquad \hbox{ when } \Phi(d(x,y))\le c_3 t.
$$
This means that both \eqref{e:hkmequi} and \eqref{eq:fibie3} enjoy
the same form
for near diagonal heat kernel estimates.
Different from \eqref{e:hkmequi} where off-diagonal heat kernel estimates are due to jump systems of the corresponding jump processes, in \eqref{eq:fibie3}  off-diagonal heat kernel estimates are characterized by diffusion property and the so-called chain arguments.
The most typical examples are the sub-Gaussian heat kernels for
Brownian motion on some fractals,
in which
$\Phi (r)=r^{d_w}$, $V(x,r)=r^d$ and $m(t,r)=(r^{d_w}/t)^{1/(d_w-1)}$
for all $x\in E$ and $t,r>0$
where $d_w\ge 2$ and $d\ge 1$.
Examples satisfying
\eqref{eq:fibie3} include diffusions on fractals such as Sierpinski gaskets and Sierpinski carpets, see
for instance \cite{bar} and \cite{HK} for details and more examples.

\begin{thm}\label{theorem:maindiff}
We assume that Assumption $\ref{as:sub}$ holds.
Suppose that the heat kernel of the conservative diffusion process
$X$ enjoys estimates \eqref{eq:fibie3}. Let $q(t,x,y)$ be given by
\eqref{e:3.3}. Then we have the following two statements:
\begin{itemize}
\item[$(i)$] If $\Phi(d(x,y))\phi(t^{-1})\le 1$, then
\begin{align*}q(t,x,y) &\simeq  \frac{ \phi(t^{-1})}{t} \int_{\Phi(d(x,y))}^{2/\phi(t^{-1})} \frac{r}{V(x,\Phi^{-1}(r))}\,dr.\end{align*}

\item[$(ii)$]If $\Phi(d(x,y))\phi(t^{-1})\ge 1$, then
there exist constants $c_i>0$ $(i=1, \dots, 4)$ such that
\begin{align*}
\frac{c_1}{tV(x, \Phi^{-1}(1/\phi(t^{-1}))))} & \frac{n(t,d(x,y))}{\phi(n(t,d(x,y))/t)}\exp(-c_2n(t,d(x,y)))\\
  & \le   q(t,x, y)\\
  &\le \frac{c_3}{tV(x, \Phi^{-1}(1/\phi(t^{-1}))))} \frac{n(t,d(x,y))}{\phi(n(t,d(x,y))/t)}\exp(-c_4n(t,d(x,y))),
\end{align*}
where
$n(\cdot, r)$ is a non-increasing function on $(0,\infty)$ determined by
\begin{align}
\label{e:phnd}
\frac{1}{\phi(n(t,r)/t)} \simeq \Phi\left(\frac{r}{n(t,r)}\right),\quad t,r>0.
\end{align}
\end{itemize}
\end{thm}

Let us roughly give some intuitive explanations on these estimates between
the fundamental solution $q(t, x, y)$ for the time fractional Poisson equation and the fundamental solution $p(r, x, y)$ of the homogenous time fractional equation.
According to \eqref{e:3.4},
$$p(t,x,y)=\int_0^\infty p^0(r,x,y)\,d_r \P(E_t\le r)=\int_0^\infty p^0(r,x,y)\,d_r\P(S_r\ge t).$$ By Proposition \ref{p:sub} and Theorem \ref{t:sub}, for any $t,r>0$,
$$ \bar{p}(r,t) \simeq t^{-1}\P(S_r\ge t).$$ Therefore,
we shall have the following rough expression
$$q(t,x,y)\approx
\int_0^\infty p^0(r,x,y)\cdot \left(\frac{r}{t}\right)\,d_r\P(S_r\ge t);$$ that is, to estimate $q(t,x,y)$ we need to add the factor $ {r}/{t}$ inside the formula above for $p(t,x,y)$. Such intuition is much clear, if we compare Theorem \ref{theorem:mainjump} (resp.\ Theorem \ref{theorem:maindiff}) with \cite[Theorem 1.6]{CKKW} (resp.\ \cite[Theorem 1.8]{CKKW}). (In particular, for the case that $\Phi(d(x,y))\phi(t^{-1})\le 1$ or the pure jump case.)

\bigskip

\noindent {\bf Proof of Theorem \ref{T:1.4}\,(ii).}\,\,
Again, by the same argument as that in  the proof of Theorem \ref{T:1.4}\,(i), we have  $\Phi(r)=r^\alpha, V(x,r)\simeq r^d$ and $\phi(r)=r^\beta$ in this case.
For the case $d(x,y)\le t^{\beta/\alpha}$ (i.e. $\Phi(d(x,y))\phi(t^{-1}) \le 1$), the
computation is the same as that of Theorem \ref{T:1.4}\,(i).
For the case $d(x,y)\ge t^{\beta/\alpha}$, according to \eqref{e:phnd}, $n=n(t,r)$ is now
determined by
\[
\frac{1}{(n/t)^\beta} \simeq (r/n)^\alpha, \quad t,r>0,
\]
hence $n\simeq (r^\alpha/t^\beta)^{1/(\alpha-\beta)}$ (note that
$0<\beta<1$ and $\alpha \geq 2$ in
this case). Since $n\ge 1$ in this case, we have $n^{1-\beta}\exp(-n)\asymp \exp(-n)$.
Given these, the conclusion is an easy consequence of Theorem \ref{theorem:maindiff}. \qed

\section{Appendix: Proofs of Theorems \ref{theorem:mainjump} and \ref{theorem:maindiff}}

In this part, we will give proofs of Theorems \ref{theorem:mainjump} and \ref{theorem:maindiff}.
Recall that ${\overline p^0}(t, x, d(x,y))$ is defined in \eqref{e:hkmequi} and \eqref{eq:fibie3}.
For simplicity, in the following we fix $x, y \in E$ and let $z=d(x, y)$. Then, ${\overline p^0}(t, x, d(x,y))={\overline p^0}(t, x, z)$.
Furthermore, we write ${\overline p^0}(t,x,z)$ and $V(x,r)$ as ${\overline p^0}(t,z)$ and $V(r)$, respectively. According to Theorem \ref{t:sub} and \eqref{e:3.3}, it holds
\begin{equation}\begin{split}\label{eellss} q(t,x,y)\asymp &\int_0^{2/\phi(t^{-1})} {\overline p^0}(r,z)\cdot \frac{1}{t}r \phi(t^{-1})\,dr\\
&+\int_{2/\phi(t^{-1})}^\infty {\overline p^0}(r,z)\cdot \frac{1}{t} \exp(-
c_0t(\phi')^{-1}(t/r))\,dr\\
=&:I_1+I_2.\end{split}\end{equation}

\noindent {\bf Proof of Theorem \ref{theorem:mainjump}.}\,\, The proof is split  into two cases.

\noindent
{\bf Case (1)}  Suppose that $\Phi(z)\phi(t^{-1})\le 1$. Then
\begin{align*}I_{1} &
\simeq
\frac{\phi(t^{-1})}{t} \left( \frac{1}{ V(z) \Phi(z)} \int_0^{\Phi(z)}r^2\,dr+  \int_{\Phi(z)}^{2/\phi(t^{-1})} \frac r{V(\Phi^{-1}(r))}\, dr \right)\\
& \simeq \frac{\phi(t^{-1})}{t}\left(\frac{\Phi(z)^2}{V(z)} +\int_{\Phi(z)}^{2/\phi(t^{-1})} \frac r{V(\Phi^{-1}(r))} \,dr\right)
\simeq\frac{\phi(t^{-1})}{t}\int_{\Phi(z)}^{2/\phi(t^{-1})}\frac r{V(\Phi^{-1}(r))} \,dr.
\end{align*}
Here in the last inequality
above we used the fact that
\begin{equation}\label{e:rrffee10}\int_{\Phi(z)}^{2/\phi(t^{-1})}\frac r{V(\Phi^{-1}(r))} \,dr\ge \int_{\Phi(z)}^{2\Phi(z)}\frac r{V(\Phi^{-1}(r))} \,dr \ge c \,\frac{\Phi(z)^2}{V(z)},\end{equation}
where the weak scaling properties of $\Phi$ and $V$ are used.

On the other hand, if  $\Phi(z)\phi(t^{-1})\le 1$, then
by changing variable $s=r \phi(t^{-1})$,
and using
 the weak scaling property of $\Phi$ with $(\alpha_1,\alpha_2)$,
 \eqref{e:phi},
 \eqref{vd1}, \eqref{e:04} and \eqref{e:03},
\begin{equation}\label{com1-100}\begin{split}
I_{2}&
\simeq \frac{1}{t\phi(t^{-1})} \int_{2}^\infty \exp(-c_0 t (\phi')^{-1}(t\phi(t^{-1})/s))  \cdot \frac 1{V(\Phi^{-1}(s/ \phi(t^{-1})))} \,ds\\
&
=\frac{1}{t\phi(t^{-1})V(\Phi^{-1}(1/\phi(t^{-1})))} \int_2^\infty \exp(-c_0 t(\phi')^{-1}(\phi'(t^{-1})/s))\cdot
\frac{V(\Phi^{-1}(1/\phi(t^{-1})))} {V(\Phi^{-1}(s/ \phi(t^{-1})) ) } \, ds\\
&\le \frac{c}{t\phi(t^{-1})V(\Phi^{-1}(1/\phi(t^{-1})))} \sum_{n=1}^\infty\int_{2^n}^{2^{n+1}} \exp(-c_0 t(\phi')^{-1}(\phi'(t^{-1})/s))\,s^{-d_1/\alpha_2}\,ds\\
&\le\frac{c}{t\phi(t^{-1})V(\Phi^{-1}(1/\phi(t^{-1})))}\sum_{n=1}^\infty \exp(-c_0 t(\phi')^{-1}(\phi'(t^{-1})/2^n))\,2^{-n((d_1/\alpha_2)-1) }\\
&\le \frac{c}{t\phi(t^{-1})V(\Phi^{-1}(1/\phi(t^{-1})))}\sum_{n=1}^\infty \exp(-c_12^{n(1-\beta_2)})\,2^{-n((d_1/\alpha_2)-1)}\\
&\le \frac{c}{t\phi(t^{-1})V(\Phi^{-1}(1/\phi(t^{-1})))}. \end{split}\end{equation}
Thus
 $$I_2\le  \frac{c}{t\phi(t^{-1})V(\Phi^{-1}(1/\phi(t^{-1})))}.$$
Note that
 \begin{equation}\label{e:rrffee1}\begin{split}
 \frac{\phi(t^{-1})}{t}
  \int_{\Phi(z)}^{2/\phi(t^{-1})}\frac r{V(\Phi^{-1}(r))} \,dr\ge&  \frac{\phi(t^{-1})}{t}\int_{1/\phi(t^{-1})}^{2/\phi(t^{-1})}\frac r{V(\Phi^{-1}(r))} \,dr \\ \ge&
  \frac{c}{t\phi(t^{-1})V(\Phi^{-1}(1/\phi(t^{-1})))}.\end{split}
 \end{equation}

 According to all the estimates above, we know that for any $z,t>0$ with $\Phi(z)\phi(t^{-1})\le 1$,
 \begin{align}
 \label{e:AA1}
 c_2\frac{\phi(t^{-1})}{t}\int_{\Phi(z)}^{2/\phi(t^{-1})}\frac r{V(\Phi^{-1}(r))} \,dr&\le I_1\le I_1+I_2\le c_3\frac{\phi(t^{-1})}{t}\int_{\Phi(z)}^{2/\phi(t^{-1})}\frac r{V(\Phi^{-1}(r))} \,dr.\end{align}

\noindent
{\bf Case (2)} Suppose that $\Phi(z)\phi(t^{-1})\ge 1$. Then
\begin{align*}I_{1}&
\simeq \frac{\phi(t^{-1})}{t}
\int_0^{2/\phi({t^{-1}})}  \frac{r^2}{V(z)\Phi(z)}\,dr
\simeq \frac{1}{t \phi({t^{-1}})^2 V(z) \Phi(z)}.
\end{align*}
On the other hand, following the same argument as
\eqref{com1-100}, we find that if $z,t>0$ with $\Phi(z)\phi(t^{-1})\ge 1,$ then
\begin{align*}I_{2}
&\le   \frac{c}{t\Phi(z)V(z)} \int_{2/\phi(t^{-1})}^\infty \exp(-c_0
t (\phi')^{-1}(t/r))\cdot r\,
dr\\
&\le  \frac{c}{t\phi(t^{-1})^2 \Phi(z)V(z)}   \int_2^\infty
\exp(-c_0 t(\phi')^{-1}(\phi'(t^{-1})/s))\,\cdot s\,ds\le \frac{c}{t\phi(t^{-1})^2 \Phi(z)V(z)}.
\end{align*}
Hence,  for any $z,t>0$ with $\Phi(z)\phi(t^{-1})\ge 1$,
\begin{align}
\label{e:AA2}
\frac{c_1}{t \phi({t^{-1}})^2 V(z) \Phi(z)}&\le I_1\le I_1+I_2\le \frac{c_2}{ t\phi({t^{-1}})^2V(z) \Phi(z)}. \end{align}

\ \

The desired estimates for $q(t,x,y)$ now follows from
\eqref{e:AA1}, \eqref{e:AA2}
and \eqref{eellss}. \qed

\noindent {\bf Proof of Theorem \ref{theorem:maindiff}.}\,\, The proof is divided into two
cases again.

\noindent
{\bf Case (1)}  Suppose that $\Phi(z)\phi(t^{-1})\le 1$.
Then by  \eqref{e:mphisim1}
\begin{align*}I_{1}&
\simeq\frac{\phi(t^{-1})}{t}\int_0^{\Phi(z)} \frac{r}{V(\Phi^{-1}(r))} \exp(-
m(r,z))\,dr+
\frac{\phi(t^{-1})}{t} \int_{\Phi(z)}^{2/\phi(t^{-1})} \frac{r}{V(\Phi^{-1}(r))}\,dr\\
&=: \frac{\phi(t^{-1})}{t}(I_{1,1}+I_{1,2}).
\end{align*} According to \eqref{vd1} and
\eqref{e:scdf-1}, we have
\begin{equation}\label{e:PF}\begin{split}
I_{1,1}=&\sum_{n=0}^\infty \int_{\Phi(z)/2^{n+1}}^{\Phi(z)/2^{n}}\frac{r}{V(\Phi^{-1}(r))} \exp(-m(r,z))\,dr\\
\le&\sum_{n=0}^\infty \frac{(\Phi(z)/2^{n})^2}{V(\Phi^{-1}(\Phi(z)/2^{n+1}))} \exp(-m(\Phi(z)/2^{n},z))\\
\le& \frac{c\Phi(z)^2}{V(z)} \sum_{n=0}^\infty 2^{n((d_2/\alpha_1)-2)} \exp(-c_1m(\Phi(z),z) 2^{n/(\alpha_2-1)}) \le \frac{c\Phi(z)^2}{V(z)},
\end{split}\end{equation}
where in the last inequality we used
\eqref{e:mphisim1}.
 This estimate along with \eqref{e:rrffee10}
yields that
$$I_1
\simeq  \frac{ \phi(t^{-1})}{t} \int_{\Phi(z)}^{2/\phi(t^{-1})} \frac r{V(\Phi^{-1}(r))} \,dr.$$

On the other hand, if $\Phi(z)\phi(t^{-1})\le 1$, then, by
the argument of \eqref{com1-100},
\begin{align*}I_{2}
&\le \frac{c }{t} \int_{2/\phi(t^{-1})}^\infty \exp(-c_0 t (\phi')^{-1}(t/r)) \cdot \frac1 {V(  \Phi^{-1}(r))} \,dr\le \frac{c}{t\phi(t^{-1})V(\Phi^{-1}(1/\phi(t^{-1})))}.
\end{align*}
We have by \eqref{e:rrffee1} that
$$I_1+I_2\le \frac{ c\phi(t^{-1})}{t} \int_{\Phi(z)}^{2/\phi(t^{-1})} \frac r{V(\Phi^{-1}(r))} \,dr .$$
Therefore, for any $z,t>0$ satisfying $\Phi(z)\phi(t^{-1})\le 1$, it holds that
\begin{align}\label{e:AA4}\frac{ c_1\phi(t^{-1})}{t} \int_{\Phi(z)}^{2/\phi(t^{-1})} \frac r{V(\Phi^{-1}(r))} \,dr&\le I_1\le I_1+I_2 \le \frac{ c_2\phi(t^{-1})}{t} \int_{\Phi(z)}^{2/\phi(t^{-1})} \frac r{V(\Phi^{-1}(r))} \,dr.\end{align}

\noindent
{\bf Case (2)} Suppose that $\Phi(z)\phi(t^{-1})\ge 1$.  Then by
\eqref{vd1} and \eqref{e:scdf-1}, we have
\begin{equation}\label{eq:oejsq}\begin{split} I_{1}\le &\frac{c\phi(t^{-1})}{t}
\int_{0}^{2/\phi(t^{-1})} \frac{r}{V(\Phi^{-1}(r))} \exp(-m(r,z))\,dr  \\
=&\frac{c\phi(t^{-1})}{t}\sum_{n=0}^\infty \int_{2^{1-n}/\phi(t^{-1})}^{2^{-n}/\phi(t^{-1})}\frac{r}{V(\Phi^{-1}(r))} \exp(-m(r,z))\,dr \\
\le& \frac{c\phi(t^{-1})}{t}\sum_{n=0}^\infty \frac{(2^{-n}/\phi(t^{-1}))^2}{V(\Phi^{-1}(2^{1-n}/\phi(t^{-1})))} \exp(-m(2^{-n}/\phi(t^{-1}),z)) \\
\le& \frac{c}{t\phi(t^{-1})V(\Phi^{-1}(1/\phi(t^{-1})))} \sum_{n=0}^\infty 2^{n((d_2/\alpha_1)-2)} \exp\left(-c_1m(1/\phi(t^{-1}),z) 2^{n/(\alpha_2-1)}\right) \\
\le& \frac{c}{t\phi(t^{-1})V(\Phi^{-1}(1/\phi(t^{-1})))}
\exp\left(-c_1m(1/\phi(t^{-1}),z)\right), \end{split}\end{equation} where in the
last inequality we used the facts that the function $r\mapsto
m(r,z)$ is
non-increasing and
\eqref{e:mphisim1} so that
\begin{align*}\exp\left(-c_1m(1/\phi(t^{-1}),z) (2^{n/(\alpha_2-1)}-1)\right)
 &\le \exp\left(-c_1m(\Phi(z),z) (2^{n/(\alpha_2-1)}-1)\right)\\
 &\le  \exp\left(-c_2 (2^{n/(\alpha_2-1)}-1)\right).\end{align*}

To get the estimate for $I_{2}$, we need to consider the following
two functions inside the estimates of $\bar{p}(r,t)$ and the
exponential terms of ${{\overline p^0}}(r,z)$ respectively:
\begin{equation}\label{o:fun}
G_1(r)=t(\phi')^{-1}(t/r) ~~\mbox{ and }~~G_2(r)=m(r,z)\end{equation}
for all $r>0$ and fixed $z,t >0$.
Note that, by
\eqref{e:03}, \eqref{e:scdf-1} and the facts that $\phi'$ and $m(\cdot,z)$ are non-increasing on $(0,\infty)$,
$G_1(r)$ is
a non-decreasing
function on $(0,\infty)$ such that $G_1(0)=0$ and $G_1(\infty)=\infty$, and $G_2(r)$ is a non-increasing function on $(0,\infty)$ such that $G_2(0)=\infty$ and
$G_2(\infty)=0$.
Thus, noting that $m(r,z)$ satisfies \eqref{e:scdf-1}, there exist $c^*,c_*>0$ and
$r_0=r_0(z,t)\in (0,\infty)$ such that
$G_1(r_0)
\simeq G_2(r_0)$, $G_1(r)\ge c^*G_2(r)$ when $r\ge r_0$, and $G_1(r)\le c_*G_2(r)$ when $r\le r_0$.

On the other hand, since $\Phi(z)\phi(t^{-1})\ge 1$,
by \eqref{e:04}, \eqref{e:mphisim1} and the fact that
$m(\cdot, z)$ is non-increasing on $(0,\infty)$,
\begin{align*}
G_1(1/\phi(t^{-1}))&=t(\phi')^{-1}(t\phi(t^{-1})) \simeq t(\phi')^{-1}(\phi'(t^{-1}))=1\\
&\le c_3m(\Phi(z),z)\le c_1m(1/\phi(t^{-1}), z)=c_1 G_2(1/\phi(t^{-1}))\end{align*}  and
\begin{align*} G_1(\Phi(z))\ge& G_1(1/\phi(t^{-1}))\simeq t(\phi')^{-1}(\phi'(t^{-1}))=1\ge c_4m(\Phi(z), z)= c_2 G_2(\Phi(z)),
\end{align*}
where constants $c_3, c_4$ are independent of $t$ and $z$.
Hence there are constants $c_5,c_6>0$ independent of $t$ and $z$ such that
\begin{equation}\label{e:ssddee} \frac{2}{\phi(t^{-1})} \le c_5 r_0\le c_6 \Phi(z).\end{equation}

Combining all the estimates above, we find that
\begin{equation*}\begin{split}I_{2}&=\frac{1}{t}\int_{2/\phi(t^{-1})}^\infty \exp(-c_0 t (\phi')^{-1}(t/r))  \cdot {{\overline p^0}}(r,z)\,dr\\
&\le
\frac{1}{ t V(\Phi^{-1}(2/\phi(t^{-1})))}\int_{2/\phi(t^{-1})}^\infty  \exp(-c_0 t (\phi')^{-1}(t/r))\cdot \exp(-m(r,z))\,dr\\
&\le  \frac{c}{tV(\Phi^{-1}(1/\phi(t^{-1})))}\int_{2/\phi(t^{-1})}^{c_5r_0}\exp(-m(r,z))\,dr\\
&\quad +\frac{c}{tV(\Phi^{-1}(1/\phi(t^{-1})))}\int_{c_5r_0}^\infty \exp(-c_0 t (\phi')^{-1}(t/r))\,dr \\
&=:\frac{c}{tV(\Phi^{-1}(1/\phi(t^{-1})))} I_{2,1}+\frac{c}{tV(\Phi^{-1}(1/\phi(t^{-1})))} I_{2,2}.   \end{split}\end{equation*}

According to \eqref{e:scdf-1},
\begin{align*}I_{2,1}&\le \int_0^{c_5r_0} \exp(-m(r,z))\,dr
=\sum_{n=0}^\infty \int_{c_5r_0/2^{n+1}}^{c_5r_0/2^n} \exp(-m(r,z))\,dr\\
&\le c r_0 \sum_{n=0}^\infty 2^{-n} \cdot \exp\left(-m(c_5r_0/2^n,z)\right)\le cr_0 \sum_{n=0}^\infty 2^{-n} \cdot \exp\left(-c_6m(r_0,z) 2^{n/(\alpha_2-1)}\right)\\
&\le  c r_0\exp(-c_6 G_2(r_0)),\end{align*} where in the last
inequality we used the facts that the function $r\mapsto m(r,z)$ is
non-increasing
 and
\begin{align*} \sum_{n=0}^\infty 2^{-n} \cdot \exp\left(-c_6m(r_0,z) (2^{n/(\alpha_2-1)}-1)\right)
& \le \sum_{n=0}^\infty 2^{-n} \cdot \exp\left(-c_7m(\Phi(z),z) (2^{n/(\alpha_2-1)}-1)\right)\\
&\le \sum_{n=0}^\infty  2^{-n} \cdot \exp\left(-c_8 (2^{n/(\alpha_2-1)}-1)\right)\le c_9, \end{align*} thanks to \eqref{e:ssddee}.
On the other hand, by \eqref{e:031},
\begin{align*}I_{2,2}&=\sum_{n=0}^\infty \int_{2^nc_5r_0}^{2^{n+1}c_5r_0} \exp\left(-c_0 t (\phi')^{-1}(t/r)\right)\,dr\\
&\le cr_0\sum_{n=0}^\infty 2^n \cdot \exp\left(-c_0 t (\phi')^{-1}(t/(2^nc_5r_0))\right)\le c r_0\sum_{n=0}^\infty 2^n \cdot\exp\left(-c_{10}t (\phi')^{-1}(t/r_0) 2^{n(1-\beta_2)}\right) \\
&\le c r_0\exp\left(-c_{10} t (\phi')^{-1}(t/r_0)\right)=c
r_0\exp\left(-c_{10}G_1(r_0)\right) \le c r_0 \exp\left(-c_{11}G_2(r_0)\right),
\end{align*} where in the last inequality we used the facts that the
function $r\mapsto m(r,z)$ is
non-increasing and
\begin{align*}  \sum_{n=0}^\infty 2^n \exp\left(-c_{10}t (\phi')^{-1}(t/r_0) (2^{n(1-\beta_2)}-1)\right)
&\le \sum_{n=0}^\infty  2^n  \exp\left(-c_{10}t (\phi')^{-1}(c_{12} t\phi(t^{-1})) (2^{n(1-\beta_2)}-1)\right)\\
&\le \sum_{n=0}^\infty  2^n  \exp\left(-c_{13}(2^{n(1-\beta_2)}-1)\right)\le c_{14}, \end{align*} thanks to \eqref{e:ssddee} again.
Putting these estimates together, we obtain
\begin{equation}\label{e:llff} \begin{split}I_{2}\le \frac{c}{tV(\Phi^{-1}(1/\phi(t^{-1})))}r_0\exp\left(-(c_6 \wedge c_{11})G_2(r_0)\right).\end{split}\end{equation}

\ \

Next, we rewrite the exponential term in the right hand side of \eqref{e:llff}. By the fact that
$m(r_0,z)=G_2(r_0)\asymp G_1(r_0)=t(\phi')^{-1}(t/r_0)$
and
\eqref{e:scdf}, we have
$$\frac{r_0}{t(\phi')^{-1}(t/r_0)}\simeq \Phi\left( \frac{z}{t(\phi')^{-1}(t/r_0)}\right).$$   Let $s_0=(\phi')^{-1}(t/r_0)$. Then $t/r_0=\phi'(s_0)$ and, by \eqref{e:04},
\begin{equation}\label{e:llfgg}\frac{1}{\phi((ts_0)/t)}= \frac{1}{\phi(s_0)}\simeq \frac{1}{\phi'(s_0)s_0}\simeq \Phi\left( \frac{z}{ts_0}\right).\end{equation}
 Thus, \eqref{e:llff} and \eqref{e:llfgg} yield that
$$ I_2 \le \frac{c}{tV(\Phi^{-1}(1/\phi(t^{-1})))} r_0\exp(-c_{15}n(t,z)),$$ where
$n(t,z):=m(r_0,z)$ satisfies
$$\frac{1}{\phi(n(t,z)/t)} \simeq \Phi\left(\frac{z}{n(t,z)}\right).$$
Since $t(\phi')^{-1}(t/r_0)
\simeq n(t,z)$, $$r_0
\simeq t/\phi'(n(t,z)/t)
\simeq n(t,z)/\phi(n(t,z)/t),$$ and so, we have  \begin{align}\label{e:AA5}I_2 \le \frac{c n(t,z)}{t\phi(n(t,z)/t)V(\Phi^{-1}(1/\phi(t^{-1})))}  \exp(-c_{15}n(t,z)). \end{align}

Following the argument above, we can also obtain that
\begin{align}\label{e:AA6}I_2\ge \frac{c n(t,z)}{t\phi(n(t,z)/t)V(\Phi^{-1}(1/\phi(t^{-1})))}  \exp(-c_{16}n(t,z)). \end{align}

Set $n:=n(t,z)$. By \eqref{e:ssddee}, we have
$n= G_2(r_0) \ge m (c_5^{-1}c_6 \Phi(z),z)
\simeq 1$. Hence,
$1/\phi(t^{-1})\le c n/\phi(n/t)$ for any $t>0$.
Using this, \eqref{eq:oejsq} and
\eqref{e:ssddee} again, we get
\begin{equation}\label{e:AA7}\begin{split}
I_1&\le \frac{c_{17}}{t\phi(t^{-1})V(\Phi^{-1}(1/\phi(t^{-1})))}
\exp\left(-c_1m(1/\phi(t^{-1}),z)\right)\\
&\le \frac{c_{18} n}{t\phi(n/t)V(\Phi^{-1}(1/\phi(t^{-1})))}
\exp(-c_{19}n).  \end{split}\end{equation}

 By \eqref{e:AA5}-- \eqref{e:AA7},
 we obtain that for any $z,t>0$ with $\Phi(z)\phi(t^{-1})\ge 1$
\begin{align}\label{e:AA8} \frac{c_{20}n}{t\phi(n/t)V(\Phi^{-1}(1/\phi(t^{-1})))}  \exp(-c_{21}n)&\le I_2\le I_1+ I_2\nn\\
&\le  \frac{c_{22}n}{t\phi(n/t)V(\Phi^{-1}(1/\phi(t^{-1})))}  \exp(-c_{23}n). \end{align}

Therefore, combining  \eqref{e:AA4}, \eqref{e:AA8} and   \eqref{eellss},
we get the desired  estimates for $p(t,x,y).$ \qed

\bigskip

\noindent \textbf{Acknowledgements.}
The research of Zhen-Qing Chen is partially supported by Simons Foundation Grant 520542 and a
 Victor Klee Faculty Fellowship at UW.
The research of Panki Kim is supported by
 the National Research Foundation of Korea (NRF) grant funded by the Korea government (MSIP)
(No.\ 2016R1E1A1A01941893).\ The research of Takashi Kumagai is supported
by JSPS KAKENHI Grant Number JP17H01093 and by the Alexander von Humboldt Foundation.\
 The research of Jian Wang is supported by the National
Natural Science Foundation of China (Nos.\ 11522106 and 11831014), the Fok Ying Tung
Education Foundation (No.\ 151002), the Program for Probability and Statistics: Theory and Application (No.\ IRTL1704) and the Program for Innovative Research Team in Science and Technology in Fujian Province University (IRTSTFJ).

The main results of this paper, including the two-sided estimates on the fundamental solution $q(t, x, y)$ for the 
time fractional Poisson equations, have been reported in several conferences including 
{\it Fractional PDEs: theory, algorithms and applications} held from June 18 to 22, 2018 at Brown University, Providence,
{\it The Fifth IMS Asia Pacific Rim Meeting} from June 26 to 29, 2018 in Singapore, 
{\it Workshop on Dynamics, Control and Numerics for Fractional PDEs} from December 5 to 7, 2018 
at San Juan, Puerto Rico,
and {\it Non Standard Diffusions in Fluids, Kinetic Equations and Probability} from December 10 to 14, 2018
at CIRM, Marseille. We thank the organizers for the invitations and participants for their interest.

\renewcommand{\baselinestretch}{1}

\vskip 0.2truein

 {\bf Zhen-Qing Chen}

Department of Mathematics, University of Washington, Seattle,
WA 98195, USA

E-mail: zqchen@uw.edu

\medskip

 {\bf Panki Kim}

Department of Mathematical Sciences and Research Institute of Mathematics,

Seoul National University,
Building 27, 1 Gwanak-ro,

Gwanak-gu,
Seoul 08826, Republic of Korea

E-mail: pkim@snu.ac.kr

\medskip

{\bf Takashi Kumagai}

Research Institute for Mathematical Sciences,
Kyoto University, Kyoto 606-8502, Japan

E-mail: kumagai@kurims.kyoto-u.ac.jp

\medskip

{\bf Jian Wang}

College of Mathematics and Informatics \& Fujian Key Laboratory of Mathematical

Analysis and Applications (FJKLMAA), Fujian Normal University, 350007 Fuzhou,
P.R. China.

E-mail: jianwang@fjnu.edu.cn

\end{document}